\documentclass[english]{article}
\usepackage[T1]{fontenc}
\usepackage[latin9]{inputenc}
\usepackage{float}
\usepackage{amsmath}
\usepackage{amssymb}
\usepackage{graphicx}
\usepackage{esint}

\makeatletter

\newcommand{\lyxdot}{.}

\newcommand{\lyxaddress}[1]{
\par {\raggedright #1
\vspace{1.4em}
\noindent\par}
}

\makeatother

\usepackage{babel}
\begin{document}

\title{WLS-ENO: Weighted-Least-Squares Based Essentially Non-Oscillatory
Schemes for Finite Volume Methods on Unstructured Meshes}

\author{Hongxu Liu, Xiangmin Jiao}

\maketitle

\lyxaddress{Department of Applied Mathematics and Statistics, Stony Brook University\\
Stony Brook, NY 11794, USA}
\begin{abstract}
ENO (Essentially Non-Oscillatory) and WENO (Weighted Essentially Non-Oscillatory)
schemes are widely used high-order schemes for solving partial differential
equations (PDEs), especially hyperbolic conservation laws with piecewise
smooth solutions. For structured meshes, these techniques can achieve
high order accuracy for smooth functions while being non-oscillatory
near discontinuities. For unstructured meshes, which are needed for
complex geometries, similar schemes are required but they are much
more challenging. We propose a new family of non-oscillatory schemes,
called \emph{WLS-ENO}, in the context of solving hyperbolic conservation
laws using finite-volume methods over unstructured meshes. WLS-ENO
is derived based on Taylor series expansion and solved using a weighted
least squares formulation. Unlike other non-oscillatory schemes, the
WLS-ENO does not require constructing sub-stencils, and hence it provides
more flexible framework and is less sensitive to mesh quality. We
present rigorous analysis of the accuracy and stability of WLS-ENO,
and present numerical results in 1-D, 2-D, and 3-D for a number of
benchmark problems, and also report some comparisons against WENO.
\end{abstract}

\section{Introduction}

Many physical phenomena, such as waves, heat conduction, electrodynamics,
elasticity, etc., can be modeled by partial differential equations.
With the development of computer technology, many numerical methods
have been designed to solve these kinds of problems over the past
decades. Among these there are finite difference methods and their
generalizations, finite volume methods, and finite element methods.

In this paper, we consider the problem of reconstructing a piecewise
smooth function, in the context of finite volume methods for hyperbolic
conservation laws. Given a geometric domain $\Omega\subseteq\mathbb{R}^{d}$,
suppose $u$ is a time-dependent piecewise smooth function over $\Omega$,
such as a density function. For any connected region $\tau$, the
$d$-dimensional conservation law can be written in the form
\begin{equation}
\int_{\tau}\frac{\partial u(\vec{x},t)}{\partial t}d\vec{x}=-\int_{\partial\tau}\vec{F}(u)\cdot d\vec{a},\label{eq:conservation_law}
\end{equation}
where $\partial\tau$ is the boundary of $\tau$, and $\vec{F}$ is
a function of $u$, corresponding to the flux. 

A finite volume method solves the problem by decomposing the domain
$\Omega$ into cells $\{\tau_{i}\mid i=1,\dots,N\}$. Let $\vert\tau_{i}\vert$
denotes the volume of $\tau_{i}$ and $\overline{u}_{i}(t)=\frac{1}{\left|\tau_{i}\right|}\int_{\tau_{i}}u(\vec{x},t)\,d\vec{x}$,
the average of $u$ over $\tau_{i}$. We obtain an equation 
\begin{equation}
\frac{d\overline{u}_{i}(t)}{\partial t}=-\left|\tau_{i}\right|\int_{\partial\tau_{i}}\vec{F}(u)\cdot d\vec{a},\label{eq:fvm}
\end{equation}
for each $\tau_{i}$. The boundary integral requires using numerical
quadrature for the flux. The integration of the flux requires reconstructing
$u$ from the cell averages $\overline{u}(t)$ in an accurate and
stable fashion, and then evaluating the reconstruction at the quadrature
points along the cell boundaries. For stability, $\vec{F}\cdot\vec{n}$
is typically replaced by a numerical flux, such as the Lax-Friedrichs
flux,
\begin{equation}
\vec{F}\cdot\vec{n}=\frac{1}{2}\left[\left(\vec{F}\left(u^{-}\right)+\vec{F}\left(u^{+}\right)\right)\cdot\vec{n}-\alpha\left(u^{+}-u^{-}\right)\right],\label{eq:flux}
\end{equation}
where $u^{-}$ and $u^{+}$ are the values of $u$ inside and outside
the cell $\tau_{i}$. The parameter $\alpha$ is a constant, and it
should be an upper bound of the eigenvalues of the Jacobian of $u$
in the normal direction.

In this context, we formulate the mathematical problem addressed in
this paper as follows: Given the cell averages $\overline{u}_{i}$
of a piecewise smooth function $u(\vec{x})$ for cell $\tau_{1},\tau_{2},\dots,\tau_{N}$,
let $h_{i}$ be some length measure of cell $\tau_{i}$. Find a polynomial
approximation $\widetilde{u}_{i}(\vec{x})$ of degree at most $p-1$
over $\tau_{i}$, such that
\begin{equation}
\Vert\widetilde{u}_{i}(\vec{x})-u_{i}(\vec{x})\Vert=\mathcal{O}(h_{i}^{p}),\hspace{1em}\vec{x}\in\tau_{i}.\label{eq:reconstruction}
\end{equation}
In other words, $\widetilde{u}_{i}(\vec{x})$ is a $p$th order accurate
approximation to $u(\vec{x})$ inside $\tau_{i}$. In the context
of hyperbolic conservation laws, $u(\vec{x})$ in (\ref{eq:reconstruction})
is equal to $u(\vec{x},t)$ in (\ref{eq:conservation_law}) at a given
$t$. For the facet between two cells, these reconstructions give
us two values $u^{-}$ and $u^{+}$, which can then be substituted
into (\ref{eq:flux}) to calculate the numerical flux. These reconstructions
must be accurate, and also must lead to stable discretizations of
the hyperbolic conservation laws when coupled with some appropriate
time integration schemes, such as TVD Runge-Kutta schemes \cite{gottlieb1998total}.

This reconstruction problem is decidedly challenging, because hyperbolic
conservation laws can produce non-smooth solutions. An approximation
scheme for smooth functions may lead to oscillations that do not diminish
as the mesh is refined, analogous to the Gibbs phenomena. Such oscillations
would undermine the convergence of the solutions. The ENO (Essentially
Non-Oscillatory) and WENO (Weighted Essentially Non-Oscillatory) schemes
\cite{ShuWeno98,JiShu96} have been successful in solving this problem.
In a nutshell, the WENO schemes use a convex combination of polynomials
constructed over some neighboring cells, with higher weights for cells
with smoother solutions and lower weights for cells near discontinuities.
As a result, these methods can achieve high-order accuracy at smooth
regions while being non-oscillatory near discontinuities. These reconstructions
can be integrated into both finite volume and finite difference methods.
With years of development, finite volume WENO schemes have been applied
to both structured and unstructured meshes and higher dimensions \cite{HuShu99WENO,levy2002fourth,Shi02TTN,tsoutsanis2011weno,lahooti2012new,titarev2004finite}.
Various attempts have been applied to improve the weights for WENO
reconstruction \cite{Henrick05MWENNO,Borges08IWENO,zhang2007new,shen2010improvement}.
Also, they have used WENO schemes in many applications, such as shock
vortex interaction \cite{Grasso00SWV}, incompressible flow problems
\cite{yang1998implicit}, Hamilton-Jacobi equations \cite{jiang2000weighted},
shallow water equations \cite{Liska99TDS}, etc.

Along the path of applying WENO schemes on unstructured meshes, tremendous
effort has been made to improve the robustness of the schemes. Early
attempts \cite{HuShu99WENO} work well for most unstructured meshes,
but some point distributions may lead to negative weights and in turn
make the schemes unstable. An extension was proposed in \cite{Shi02TTN}
to mitigate the issue, but it still had limited success over complicated
geometries due to inevitably large condition numbers of their local
linear systems. More recently, several different partition techniques
were proposed to improve stability, such as \cite{liuzhang2013robust},
which uses a hybrid of two different reconstruction strategies to
achieve better results. The technique was adopted in \cite{Christlieb2015MPP,CharestENOvertex,ZhuHJunstr}
for further development.

In this paper, we propose a new family of reconstruction methods over
unstructured meshes. We refer to the schemes as \emph{WLS-ENO}, or
\emph{Weighted-Least-Squares based Essentially Non-Oscillatory }schemes.
Unlike the WENO scheme, our approach uses a generalized finite difference
(GFD) formulation based on weighted least squares, rather than a weighted
averaging of traditional finite differences. The GFD method is derived
rigorously from Taylor series, and hence can deliver the same order
of accuracy as traditional finite differences. In WLS-ENO, the convexity
requirement is satisfied automatically, since the weights are specified
\emph{a priori}. These properties enable a more systematic way to
construct non-oscillatory schemes. We will present the detailed derivation
of the schemes and their robust numerical solution techniques. We
also show that the schemes are often more accurate than WENO schemes
near discontinuities and enable more stable PDE solvers when used
in conjunction with total variation-diminishing time-integration schemes
such as TVD Runge-Kutta. We report theoretical analysis in 1-D as
well as experimental results in 1-D, 2-D, and 3-D.

The remainder of this paper is organized as follows. Section~\ref{sec:background}
reviews the ENO and WENO schemes, as well as some related background
knowledge. Section~\ref{sec:WLS-ENO} presents the derivation and
numerical methods of the WLS-ENO schemes. Section~\ref{sec:Analysis-of-LS-WENO}
analyzes the accuracy and stability of the WLS-ENO schemes, and compares
them against WENO and its previous generalization to unstructured
meshes. Section~\ref{sec:Numerical-Results} presents some numerical
results and comparisons against some other methods. Finally, Section~\ref{sec:Conclusions-and-Future}
concludes the paper with some discussions on future research directions.

\section{\label{sec:background}Background and Related Work}

In this section, we present some background information, including
WENO schemes and some related numerical methods for hyperbolic conservation
laws. These will motivate the derivation of WLS-ENO.

Traditionally, two approaches have been used to reduce oscillations
near discontinuities. One approach was to add artificial viscosity
\cite{Sod85,LeV92}.  Specifically, one could design the viscosity
to be larger near discontinuity to suppress oscillations and to be
smaller elsewhere to maintain accuracy. However, the parameter controlling
the artificial viscosity is very problem dependent. Another approach
was to apply limiters \cite{Sod85,LeV92},  but such schemes degenerate
to first order near discontinuities. Other successful methods include
ENO scheme \cite{HaEnOsCha87} and its closely related WENO schemes
\cite{ShuWeno98}. In a nutshell, the ENO is a WENO scheme with only
zeros and ones as the weights. In the following, we review some of
the ENO and WENO schemes that are most closely related to our proposed
approach.

\subsection{WENO Reconstructions in 1-D}

In the context of finite volume methods, the basic idea of WENO is
to first construct several stencils for each cell and local polynomials
over these stencils, so that the cell averages of these polynomials
are the same as the given values. Then, a WENO scheme uses a convex
combination of these polynomials to obtain a reconstruction of the
function, where the weights for each stencil are controlled by a smoothness
indicator. We briefly describe the WENO scheme on a uniform 1-D grid
below, and refer readers to \cite{ShuWeno98} for more detail.

Given a 1-D domain $[a,b]$, suppose we have a uniform grid with nodes
\begin{equation}
a=x_{\frac{1}{2}}<x_{\frac{3}{2}}<x_{\frac{5}{2}}<\dots<x_{N-\frac{1}{2}}<x_{N+\frac{1}{2}}=b.
\end{equation}
We denote $i$th cell $\left[x_{i-\frac{1}{2}},x_{i+\frac{1}{2}}\right]$
as $\tau_{i}$ for $i=1,2,\dots,N$. Its cell center is $x_{i}=\frac{1}{2}\left(x_{i-\frac{1}{2}}+x_{i+\frac{1}{2}}\right),$
and its cell size is $h_{i}=x_{i+\frac{1}{2}}-x_{i-\frac{1}{2}}$.
The cell average of a function $u(x)$ over $\tau_{i}$ is then
\begin{equation}
\overline{u}_{i}=\frac{1}{\Delta x_{i}}\int_{x_{i-\frac{1}{2}}}^{x_{i+\frac{1}{2}}}u(x)\,dx,\hspace{1em}i=1,2,\dots,N.
\end{equation}
For each cell $\tau_{i}$, our goal is to reconstruct a piecewise
polynomial approximation $\widetilde{u}_{i}(x)$ of degree at most
$p-1$, such that it approximates $u(x)$ to $p$th order accuracy
within $\tau_{i}$, i.e.,
\begin{equation}
\widetilde{u}_{i}(x)=u(x)+\mathcal{O}(h^{p}),\hspace{1em}x\in\tau_{i},\hspace{1em}i=1,\dots,N,\label{eq:reconstruction-1D}
\end{equation}
where $h=\min\{h_{i}\mid1\leq i\leq N\}$.

To find such a polynomial, a WENO scheme first selects $p$ sub-stencils
about $\tau_{i}$, each containing $p$ cells. Consider a particular
sub-stencil
\begin{equation}
S_{j}(i)=\left\{ \tau_{i-j},\dots,\tau_{i-j+p-1}\right\} ,
\end{equation}
and let $\phi_{i,j}(x)$ be a polynomial approximation of $u$ over
$S_{j}(i)$, obtained by requiring the integral of $\phi_{i,j}(x)$
over each cell in the sub-stencil to be equal to that of $u_{i}(x)$.
If the $p$th derivative of $u$ is bounded over the sub-stencil $S_{j}(i)$,
then $\phi_{i,j}(x)$ satisfies (\ref{eq:reconstruction-1D}). However,
if $u(x)$ has discontinuities within the sub-stencil, then $\phi_{i,j}(x)$
may be oscillatory. The WENO scheme then constructs a non-oscillatory
approximation by taking a convex combination of $\phi_{i,j}(x)$
\begin{equation}
\widetilde{u}_{i}(x)=\sum_{j}\omega_{j}\phi_{i,j}(x),
\end{equation}
where $\omega_{j}=\alpha_{j}/\sum_{k=0}^{p-1}\alpha_{k}$ and is chosen
such that $\omega_{j}$ approaches zero for sub-stencils with discontinuities.
A typical choice of $\alpha_{j}$ is
\begin{equation}
\alpha_{j}=d_{j}/\left(\epsilon+\beta_{j}\right)^{2},\label{eq:base_weight}
\end{equation}
where $d_{j}$ is a nonnegative coefficient such that
\begin{equation}
\widetilde{u}_{i}(x)=\sum_{j=0}^{p-1}d_{j}\phi_{i,j}(x)=u\left(x_{i+\frac{1}{2}}\right)+\mathcal{O}(h^{2p-1}).\label{eq:WENO linear scheme}
\end{equation}
The parameter $\epsilon$ is a small parameter, such as $\epsilon=10^{-6}$,
introduced to avoid instability due to division by zero or too small
a number. The non-negativity of $d_{i}$ is important for stability
purposes. The $\beta_{j}$ is the \emph{smoothness indicator}. If
$u(x)$ is smooth over the sub-stencil $S_{j}(i)$, then $\beta_{j}=\mathcal{O}(h^{2})$;
otherwise, $\beta_{j}=\mathcal{O}(1)$. A typical choice of $\beta_{j}$,
as introduced in \cite{ShuWeno98}, is
\begin{equation}
\beta_{j}=\sum_{k=1}^{p-1}\int_{x_{i-\frac{1}{2}}}^{x_{i+\frac{1}{2}}}h^{2k-1}\left(\frac{\partial^{k}\phi_{j}(x)}{\partial x^{k}}\right)^{2}\,dx,
\end{equation}
where the $h^{2k-1}$ term is introduced to make $\beta_{j}$ independent
of the grid resolution. For example, in the simplest case where $p=2$,
\begin{equation}
\begin{aligned}\beta_{0} & =\left(\overline{u}_{i+1}-\overline{u}_{i}\right)^{2}\\
\beta_{1} & =\left(\overline{u}_{i}-\overline{u}_{i-1}\right)^{2}
\end{aligned}
.
\end{equation}
Alternative smoothness indicators have been proposed in \cite{Henrick05MWENNO,Borges08IWENO,Castro11HOWENO,ZhangShu07NSI,Shen10IWS}.

\subsection{WENO Schemes on 2-D and 3-D Structured Meshes}

Originally developed in 1-D, the WENO schemes can be generalized to
structured meshes in 2-D and 3-D. Here, we give a brief overview of
the reconstructions in 2-D, which generalize to 3-D in a relatively
straightforward manner.

Consider an $N_{x}$-by-$N_{y}$ structured grid, and let $\tau_{ij}$
denotes the cell $(i,j)$ in the grid, where $i=1,2,\dots,N_{x}$
and $j=1,2,..,N_{y}$. Suppose the cell averages of a function $u(x,y)$,
\begin{equation}
\overline{u}_{ij}=\frac{1}{\Delta x_{i}\Delta y_{j}}\int_{y_{j-\frac{1}{2}}}^{y_{j+\frac{1}{2}}}\int_{x_{i-\frac{1}{2}}}^{x_{i+\frac{1}{2}}}u\left(x,y\right)\,dxdy,
\end{equation}
are given. We would like to find
\begin{equation}
\widetilde{u}_{ij}(x,y)=\sum_{r=0}^{p-1}\sum_{s=0}^{p-1}a_{rs}x^{r}y^{s},
\end{equation}
where $a_{rs}$ are the coefficients to be determined so that $\widetilde{u}_{ij}$
approximates $u$ to $p$th order accuracy over $\tau_{ij}$. 

Similar to the 1-D case, a WENO scheme first constructs polynomial
approximations over a selection of sub-stencils and then computes
a convex combination of these approximations. Consider a particular
sub-stencil
\begin{equation}
S_{lm}\left(i,j\right)=\left\{ \tau_{IJ}:i-l\leq I\leq i-l+p-1,j-m\leq J\leq j-m+p-1\right\} .
\end{equation}
Let $\phi_{ijlm}(x,y)$ denotes a polynomial reconstruction of $u$
over $S_{lm}\left(i,j\right)$, whose integration over each cell in
$S_{lm}\left(i,j\right)$ is equal to the given cell average. The
function $\phi_{ijlm}(x,y)$ approximates $u(x,y)$ over $\tau_{ij}$
to $p$th order accuracy for smooth functions. To obtain a non-oscillatory
reconstruction for non-smooth functions, the WENO scheme computes
$\widetilde{u}_{ij}$ as
\begin{equation}
\widetilde{u}_{ij}(x,y)=\sum_{l}\sum_{m}\omega_{lm}\phi_{ijlm}(x,y).
\end{equation}
The weights $\omega_{lm}$ are chosen so that the order of accuracy
of $\widetilde{u}_{ij}$ is maximized for smooth functions, and then
further augmented based on similar smoothness indicators as in 1-D,
so that the weights would approach zero for the sub-stencils with
discontinuities. For stability, it is important that the weights are
nonnegative, which imposes some constraints to the selection of stencils.
For more detail, see \cite{ShuWeno98,wolf2007high,levy2002fourth}.

\subsection{ENO and WENO Schemes on Unstructured Meshes}

Besides structured meshes, WENO can also be generalized to unstructured
meshes. Typically, these schemes are also constructed from some convex
combination of lower-order schemes on sub-stencils. However, compared
to structured meshes, it is much more challenging to construct stable
WENO schemes on unstructured meshes, because it is difficult to satisfy
the convexity requirement. Some WENO schemes have been proposed for
2-D \cite{HuShu99WENO,Shi02TTN,wolf2007high,levy2002fourth} and 3-D
\cite{yang1998implicit,tsoutsanis2011weno,zhang2009third,lahooti2012new,titarev2004finite}.
Below, we briefly review these generalizations, focusing on three
different types. 

The first type of WENO reconstruction, as proposed in \cite{HuShu99WENO},
uses a combination of high-order polynomials computed from low-order
polynomials over sub-stencils. For example, to obtain a third-order
reconstruction, the scheme would first construct linear approximations
over several sub-stencils by requiring the cell averages of the polynomials
to be equal to the given cell averages, and then compute a quadratic
polynomial from a weighted average of these linear polynomials. This
technique works for unstructured meshes, even for meshes with mixed
types of elements. However, depending on the mesh, the weighted average
may not form a convex combination, and the linear system for calculating
the weights may be ill-conditioned.

The second type of WENO reconstruction is similar to the first type,
except that it compromises the order of accuracy of the convex combination,
by allowing convex combination to be the same degree polynomials as
those for the sub-stencils. Compared to the first type, this construction
is less sensitive to mesh quality than the first type at the cost
of lower accuracy. Therefore, it is often used as a fallback of the
first type for robustness \cite{liuzhang2013robust}.

The third type of WENO reconstruction builds the reconstructions in
a hierarchical fashion. An example is the approach in \cite{wolf2007high},
which first finds the smoothest linear reconstructions over the first-layer
three-cell stencils, and then use these linear reconstructions to
build quadratic reconstructions over the second-layer stencils. This
approach can be applied iteratively to construct higher-order reconstructions.
However, as shown in \cite{wolf2007high}, the accuracy of the reconstruction
may not improve as the degree of the polynomial increases, especially
near boundaries.

Besides the above WENO schemes, we also note some recent development
of the Central ENO (CENO) schemes \cite{CharestCENO,CharestENOvertex},
based on least squares approximations. However, CENO schemes require
limiters for linear reconstructions to preserve monotonicity. Our
proposed WLS-ENO scheme differs from the WENO and CENO reconstructions
in that it utilizes a weighed least squares formulation, it does not
require limiters, and it is insensitive to mesh quality due to adaptive
stencils.

\section{Weighted-Least-Squares Based ENO Schemes \label{sec:WLS-ENO}}

In this section, we propose a new class of essentially non-oscillatory
schemes, referred to as \emph{WLS-ENO}. In the context of finite volume
methods, these schemes reconstruct a function $u(\vec{x})$ over each
cell, given the cell averages of $u$, denoted by $\overline{u}$,
for all the cells. For each point along cell boundaries, the reconstructions
then provide two values, $u^{-}$ and $u^{+}$, which can then be
used in (\ref{eq:flux}) to calculate fluxes. Unlike WENO schemes,
the WLS-ENO does not use weighted averaging of functional approximations
over sub-stencils, but instead computes the reconstruction over each
cell based on weighted least squares with an adaptive stencil. It
can achieve optimal accuracy for smooth functions, stability around
discontinuities, and insensitivity to mesh quality. We will present
the derivations of WLS-ENO schemes based on Taylor series expansion,
as well as their implementations in 1-D, 2-D and 3-D over structured
and unstructured meshes.

\subsection{WLS-ENO Schemes in 1-D}

We first derive the WLS-ENO schemes in 1-D. Suppose we are given a
grid
\begin{equation}
a=x_{\frac{1}{2}}<x_{\frac{3}{2}}<x_{\frac{5}{2}}<\cdots<x_{N-\frac{1}{2}}<x_{N+\frac{1}{2}}=b,
\end{equation}
and the cell averages $\overline{u}_{i}$ of a function $u(x)$ over
each cell $\tau_{i}$, $i=1,2,\dots,N$. For generality, we assume
the grid is non-uniform, with varying cell sizes. Below, we first
describe how to reconstruct $u$ from $\overline{u}_{i}$ for smooth
functions, and then augment the method to handle discontinuities.

\subsubsection{WLS-Based Reconstruction for Smooth Functions}

Without a loss of generality, let us consider the reconstruction of
$u$ over $\tau_{i}$ at its boundary point $x_{i+\frac{1}{2}}$.
To achieve $p$th order accuracy, we need to construct a polynomial
approximation of degree at least $p-1$. We choose a stencil with
$n$ cells to perform the reconstruction, where $n\geq p$. 

Suppose there are $l$ cells to the left of $\tau_{i}$ in the stencil.
The full stencil is given by the set
\begin{equation}
S(i)=\left\{ \tau_{i-l},\dots,\tau_{i-l+n-1}\right\} .
\end{equation}
From Taylor series expansion, we can approximate function $u(x)$
at point $x_{i+\frac{1}{2}}$ to $p$th order accuracy by
\begin{equation}
u(x)=\sum_{k=0}^{p-1}\frac{u^{(k)}(x_{i+\frac{1}{2}})}{k!}\left(x-x_{i+\frac{1}{2}}\right)^{k}+\mathcal{O}(h^{p}),
\end{equation}
where $h$ denotes the average edge lengths. The cell average over
$\tau_{j}$ in the stencil can be approximated by
\begin{align*}
\overline{u}_{j} & =\sum_{k=0}^{p-1}\frac{u^{(k)}(x_{i+\frac{1}{2}})}{k!\left(x_{j+\frac{1}{2}}-x_{j-\frac{1}{2}}\right)}\int_{x_{j-\frac{1}{2}}}^{x_{j+\frac{1}{2}}}\left(x-x_{i+\frac{1}{2}}\right)^{k}\,dx+\mathcal{O}(h^{p})\\
 & =\sum_{k=0}^{p-1}\frac{u^{(k)}(x_{i+\frac{1}{2}})}{(k+1)!}\left[\left(x_{j+\frac{1}{2}}-x_{i+\frac{1}{2}}\right)^{k+1}-\left(x_{j-\frac{1}{2}}-x_{i+\frac{1}{2}}\right)^{k+1}\right]+\mathcal{O}(h^{p}).
\end{align*}
Given $\overline{u}_{j}$ at the $j$th cells in $S(i)$, we then
construct an $n\times p$ linear system
\begin{equation}
\vec{A}\vec{v}\approx\overline{\vec{u}},\label{eq:average_eq}
\end{equation}
where 
\begin{equation}
a_{JK}=\frac{1}{K!}\left[\left(x_{i-l+J-\frac{1}{2}}-x_{i+\frac{1}{2}}\right)^{K}-\left(x_{i-l+J-\frac{3}{2}}-x_{i+\frac{1}{2}}\right)^{K}\right]
\end{equation}
for $J\equiv j+l+1-i\in[1,n]$ and $K\equiv k+1\in[1,p]$, $\overline{\vec{u}}$
is composed of the cell averages $\overline{u}_{j}$, and $\vec{v}$
is composed of the derivative of function $u$ at $x_{i+\frac{1}{2}}$,
i.e., $v_{K}=u^{(k)}(x_{i+\frac{1}{2}})$.

Eq.~(\ref{eq:average_eq}) in general is a rectangular linear system,
and we can solve it using a weighted least squares formulation. In
particular, we assign a different weight to each cell. Let $\vec{W}$
denotes a diagonal matrix containing these weights. The problem can
be written in matrix form as
\begin{equation}
\min\Vert\vec{W}\vec{A}\vec{v}-\vec{W}\overline{\vec{u}}\Vert_{2}.\label{eq:WLS}
\end{equation}
The weights allow us to assign different priorities to different cells.
For example, we may give higher weights to the cells closer to $I_{i}$.
We solve this weighted least squares problem using QR factorization
with column pivoting, as we will describe in more detail in Section~\ref{sub:Implementation-Details}.
Since the method is derived based on Taylor series expansions directly,
this WLS-based reconstruction can deliver the same order of accuracy
as interpolation-based schemes for smooth functions, as we demonstrate
in Section~\ref{sub:Accuracy}. For discontinuous functions, these
weights can also allow us to suppress the influence of cells close
to discontinuities, as we discuss next.

\subsubsection{WLS-ENO for Discontinuous Functions}

To apply WLS-ENO schemes to discontinuous functions, we modify the
weighting matrix $\vec{W}$ in (\ref{eq:WLS}). The main idea is to
assess the smoothness of the function within each cell of the stencil,
and then define the weights correspondingly. By letting the weights
be far smaller for the cells near discontinuities than those away
from discontinuities, we can then effectively suppress oscillations.

We first construct a \emph{non-smoothness} \emph{indicator} of the
function, analogous to those used in WENO schemes. Specifically, for
the $j$th cell in the stencil for $I_{i}$, with $j=i-l+J$, we can
define the indicator for cell $I_{j}$ as
\begin{equation}
\beta_{j}=\begin{cases}
\left(\overline{u}_{j}-\overline{u}_{i}\right)^{2}+\epsilon h^{2} & j\neq i\\
\min\left\{ \beta_{j-1},\beta_{j+1}\right\}  & j=i
\end{cases},\label{eq:smoothness indicator}
\end{equation}
where $\epsilon$ is a small constant, such as $\epsilon=10^{-2}$,
introduced to avoid the indicator being too close to zero, and $h$
is some measure of average edge length. Note that $\beta_{j}=\mathcal{O}(h^{2})$
if $u$ is smooth at $\tau_{j}$ and $\beta_{j}=\mathcal{O}(1)$ near
discontinuities. Therefore, it captures the non-smoothness of the
function. We therefore refer to $\beta$ as a non-smoothness indicator,
although its counterpart in WENO is called the ``smoothness indicator.''

We then define the weights based on $\beta_{j}$. To suppress oscillations,
it is desirable to use smaller weights for cells at discontinuities.
Therefore, we make the weights in $\vec{W}$ to be inversely proportional
to $\beta_{j}$ when $j\neq i$, and make the value larger for $\beta_{i}$.
Specifically, we choose
\begin{equation}
\omega_{j}=\begin{cases}
1/\beta_{j} & j\neq i\\
\alpha/\beta_{j} & j=i
\end{cases},\label{eq:weight}
\end{equation}
where $\alpha>1,$ such as $\alpha=1.5$. It is easy to see that $\omega_{j}$
is $\mathcal{O}(1/h^{2})$ if the function is smooth around cell $i$
and $\omega_{q}=\mathcal{O}(1)$ if the function is discontinuous
in the cell. After computing the weights, we substitute them into
(\ref{eq:WLS}) to compute the reconstruction. Note that unlike the
weights in WENO, we do not need to normalize the weights by dividing
them by the sum of the weights. As we will demonstrate in Section~\ref{sub:Accuracy},
this approach effectively suppresses the oscillations near discontinuities,
similar to WENO schemes.

\subsection{Generalization of WLS-ENO Schemes to 2-D and 3-D}

The WLS-ENO reconstruction can be generalized to 2-D and 3-D, for
evaluating the values at quadrature points along the cell boundaries.
Similar to 1-D, we derive the higher-dimensional version of the linear
system (\ref{eq:average_eq}) over each cell based on Taylor series
expansion, and then solve it based on weighted least squares.

Let us first consider the scheme in 2-D. First, we choose $n$ cells
with index $i_{1},i_{2},\dots,i_{n}$ as the stencil for cell $\tau_{i}$.
Let $\left(x_{i},y_{i}\right)$ denotes its centroid. From the 2-D
Taylor series expansion, we can approximate $u(x,y)$ about $\left(x_{i},y_{i}\right)$
to $p$th order accuracy by 
\begin{equation}
u(x,y)=\sum_{q=0}^{p-1}\sum_{k,l\geq0}^{k+l=q}\frac{\partial^{q}u\left(x_{i},y_{i}\right)}{\partial x^{k}\partial y^{l}}\frac{\left(x-x_{i}\right)^{k}\left(y-y_{i}\right)^{l}}{k!l!}+O\left(\left\Vert \delta\right\Vert ^{p}\right),\label{eq:2D taylor}
\end{equation}
where $\delta=\max\left\{ \left|x-x_{i}\right|,\left|y-y_{i}\right|\right\} .$
Thus, the cell averages over the $j$th cell $\tau_{j}$ can be approximated
by 
\begin{equation}
\overline{u}_{j}=\frac{1}{k!l!\left|\tau_{j}\right|}\sum_{q=0}^{p-1}\sum_{k,l\geq0}^{k+l=q}\frac{\partial^{q}u\left(x_{i},y_{i}\right)}{\partial x^{k}\partial y^{l}}\iint_{\tau_{j}}\left(x-x_{i}\right)^{k}\left(y-y_{i}\right)^{l}\,dxdy+O\left(\left\Vert \delta\right\Vert ^{p}\right)
\end{equation}
for $j=i_{1},\dots,i_{n}$, where $\left|\tau_{j}\right|$ denotes
the area of $\tau_{j}$. Therefore, we obtain $n$ equations from
the $n$ cells in the stencil about $\tau_{i}$, which can then be
solved using the weighted least squares formulation (\ref{eq:WLS}).

The construction in 3-D is based on the 3-D Taylor series expansion
about a centroid $(x_{i},y_{i},z_{i})$ of the $i$th cell $\tau_{i}$,
\begin{equation}
u(x,y,z)=\sum_{q=0}^{p-1}\sum_{k,l,m\geq0}^{k+l+m=q}\frac{\partial^{q}u\left(x_{i},y_{i},z_{i}\right)}{\partial x^{k}\partial y^{l}\partial z^{m}}\frac{\left(x-x_{i}\right)^{k}\left(y-y_{i}\right)^{l}\left(z-z_{i}\right)^{m}}{k!l!m!}+O\left(\left\Vert \delta\right\Vert ^{p}\right),\label{eq:3D taylor-1}
\end{equation}
where $\delta=\max\left\{ \left|x-x_{i}\right|,\left|y-y_{i}\right|,\left|z-z_{i}\right|\right\} .$
Then the cell averages over the $j$th cell $\tau_{j}$ in the stencil
for $\tau_{i}$ can be approximated by
\begin{equation}
\overline{u}_{j}=\frac{1}{k!l!m!\left|\tau_{j}\right|}\sum_{q=0}^{p-1}\sum_{k,l,m\geq0}^{k+l+m=q}\frac{\partial^{q}u\left(x_{i},y_{i},z_{i}\right)}{\partial x^{k}\partial y^{l}\partial z^{m}}\iiint_{\tau_{j}}\left(x-x_{i}\right)^{k}\left(y-y_{i}\right)^{l}\left(z-z_{i}\right)^{m}\,dxdydz+O\left(\left\Vert \delta\right\Vert ^{p}\right)
\end{equation}
for $j=i_{1},\dots,i_{n}$, where $\left|\tau_{j}\right|$ denotes
the volume of $\tau_{i}$. We can then solve the resulting least squares
problem using weighted least squares.

To determine the weights in (\ref{eq:WLS}), we define the non-smoothness
indicator $\beta_{j}$ similar to (\ref{eq:smoothness indicator}),
and then define the weights $\omega_{j}$ as in (\ref{eq:weight}).
As in 1-D, it is easy to show that $\beta_{j}=\mathcal{O}(h^{2})$
if $u$ is smooth at $\tau_{j}$ and $\beta_{j}=\mathcal{O}(1)$ near
discontinuities, and therefore the weights can effectively suppress
the effect of cells near discontinuities. After solving the linear
system and obtaining the polynomial approximation, we can then evaluate
$u(x,y)$ at the quadrature points along the boundary for computing
the numerical fluxes. 

One remaining question is the selection of the stencils. For the reconstruction
over triangular meshes in 2-D, we adopt the strategy in \cite{JiaoWang10RHO}
to define $k$-ring neighbor cells, with $1/2$-ring increments:
\begin{itemize}
\item The $1/2$-ring neighbor cells are the cells that share an edge with
a cell.
\item The 1-ring neighbor cells of a cell are those that share at least
one vertex with the center cell.
\item For any positive integer $k$, the $(k+1)$-ring neighborhood of a
cell is the union of $k$-ring neighborhood and 1-ring neighbors of
its $k$-ring neighbor cells. The $(k+1/2)$-ring neighborhood is
the union of $k$-ring neighborhood and $1/2$-ring neighbors of the
$k$-ring neighbor cells.
\end{itemize}
Figure~\ref{fig:Examples-of-1-ring,} illustrates the neighborhood
definitions up to 2 rings. The $1/2$-ring increment allows finer
granularity in the increment of the stencil sizes. 
\begin{figure}
\begin{centering}
\includegraphics[scale=0.4]{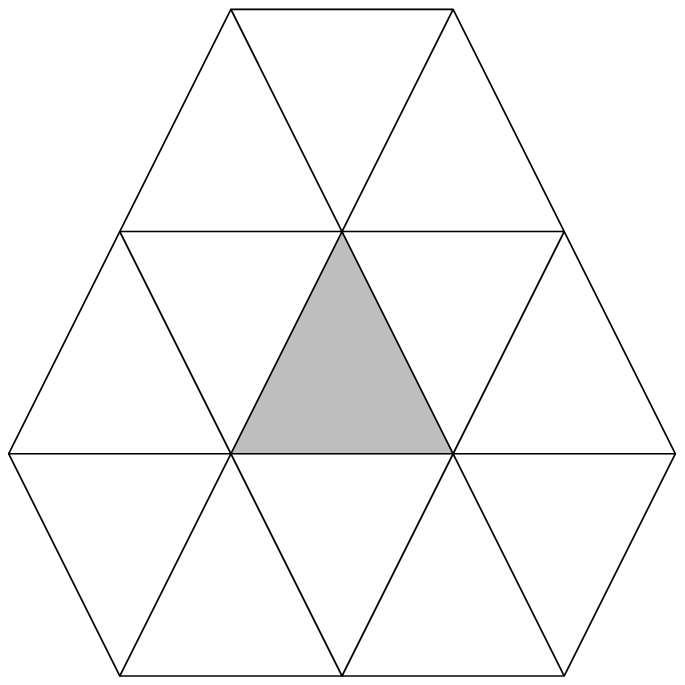} \includegraphics[scale=0.4]{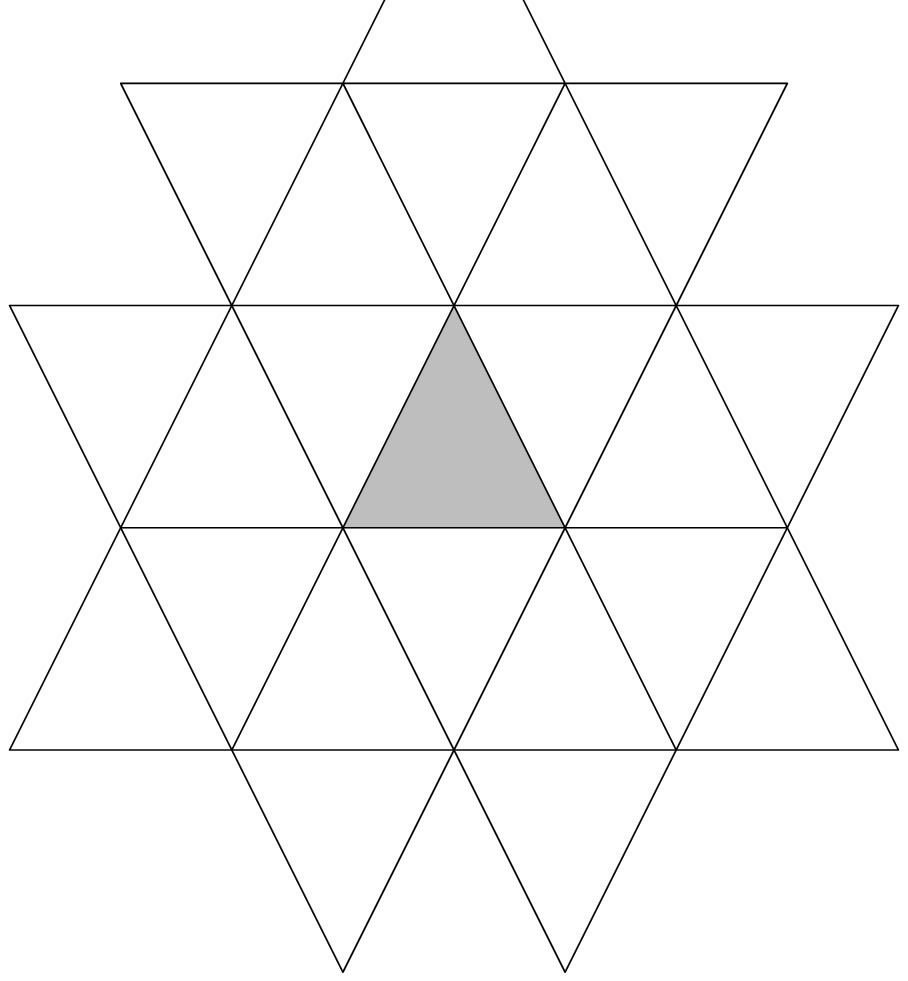}\includegraphics[scale=0.4]{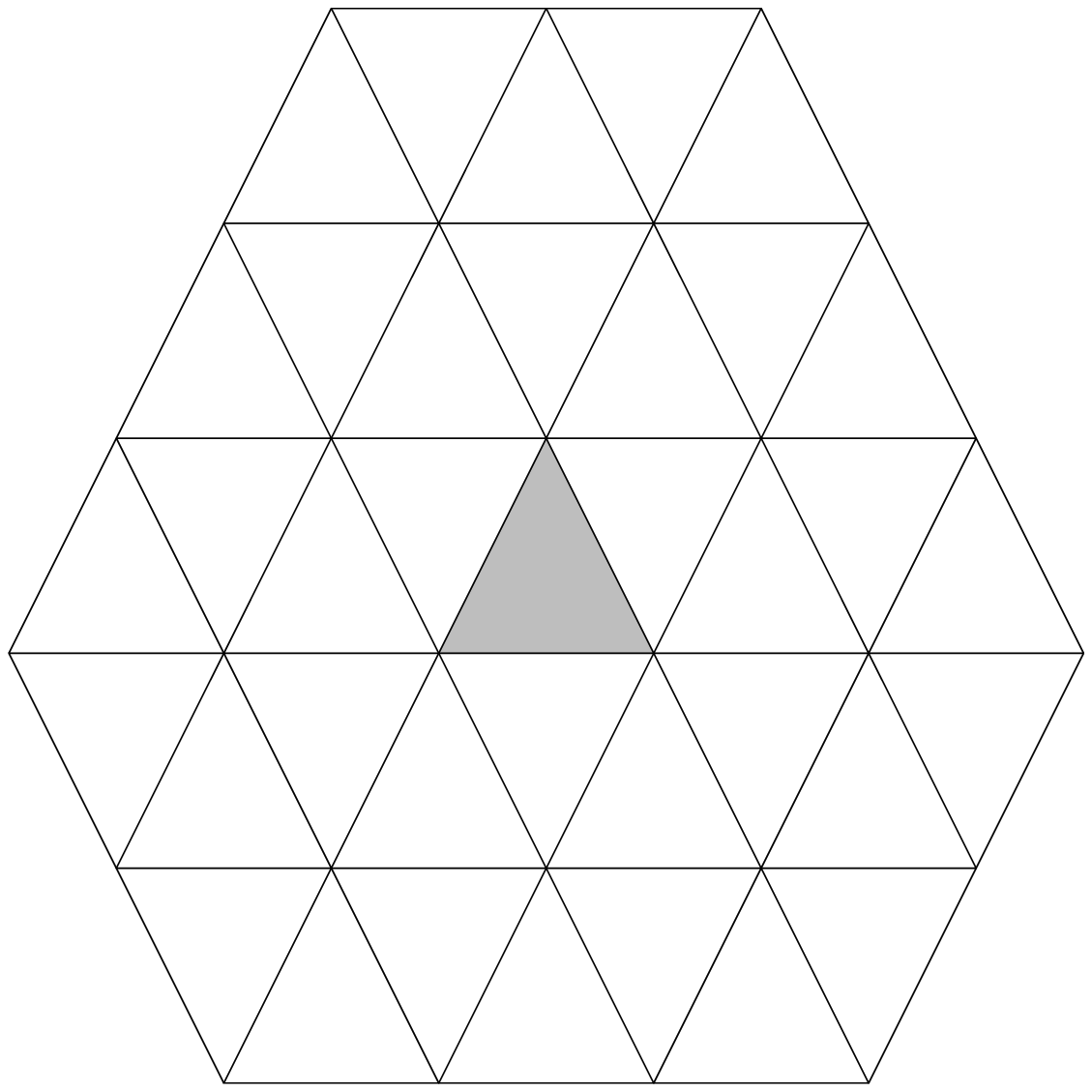}
\par\end{centering}

\protect\caption{\label{fig:Examples-of-1-ring,}Examples of 1-ring, 1.5-ring and 2-ring
neighborhood of a triangle.}
\end{figure}

For the reconstruction over tetrahedral meshes in 3-D, the standard
$k$-ring neighbors grow very rapidly. To allow finer granularity,
we define $k$-ring neighbor cells with $1/3$-ring increments, similar
to those defined in \cite{Conley16AESFEM}:
\begin{itemize}
\item The $1/3$-ring neighbor cells of a cell are the cells that share
at least one face with the center cell.
\item The $2/3$-ring neighbor cells of a cell are the cells that share
at least one edge with the center cell.
\item The 1-ring neighbor cells of a cell are the cells that share at least
one vertex with the center cell.
\item For any positive integer $k$, the $\left(k+1/3\right)$-ring neighborhood
of a cell is the union of the $k$-ring neighborhood and the $1/3$-ring
neighbors of its $k$-ring neighbor cells. The $\left(k+2/3\right)$-ring
neighborhood is the union of the $k$-ring neighborhood and the $2/3$-ring
neighbors of the $k$-ring neighbor cells.
\end{itemize}
With the above definitions, we adaptively choose the $k$-ring neighborhood,
so that the number of cells in a stencil is approximately equal to
1.5 to 2 times the number of coefficients in the Taylor polynomial.
This adaptive strategy allows the WLS-ENO scheme to be less dependent
on mesh quality than WENO schemes.

The above strategy works well in practice for most cases. However,
on poor-quality meshes, some $k$-ring neighborhoods may be nearly
one-sided, which may cause the least-squares approximation to be closer
to extrapolation along some directions. This may cause the reconstructed
value to fall beyond the maximum and minimum values of the cell averages
in the neighborhood, and in turn may lead to oscillations. This rarely
happens in 1-D or 2-D, but we do observe it in practice for some 3-D
meshes. This issue could be mitigated by using limiters, analogous
to the approach in \cite{Parklimiting} for multidimensional reconstructions
on unstructured meshes. However, we resolve the issue by adapting
the stencil as follows. First, after the reconstruction, we check
whether the reconstructed values are between the maximum and minimum
values in the neighborhood. If not, we compute the plane that passes
through the centroid $\vec{x}_{i}$ of $\tau_{i}$ and is orthogonal
to $\nabla u$ at $\vec{x}_{i}$.  Next, we select a subset of the
cells from an enlarged stencil to ensure the new stencil is well balanced
on the two sides of the plane. Because of the smoothing nature of
least squares, we find that the new polynomial approximation typically
falls within the range on balanced stencils.

\subsection{\label{sub:Implementation-Details}Implementation Details}

The WLS-ENO is different from the WENO schemes in terms of the stencil
selection and the local linear systems, and therefore it requires
different data structures and linear algebra techniques. Below, we
address some of these implementation details.

\subsubsection{Data Structure for Neighborhood Search}

To support the construction of stencils in 2-D and 3-D, we use an
Array-based Half-Facet (AHF) data structure \cite{DREJTAHF2014} to
store the mesh information. In a $d$-dimensional mesh, the term \emph{facet}
refers to the $(d-1)$-dimensional mesh entities; that is, in 2-D
the facets are the edges, and in 3-D the facets are the faces. The
basis for the half-facet data structure is the idea that every facet
in a manifold mesh is made of two half-facets oriented in opposite
directions. We refer to these two half-facets as \emph{sibling half-facets}.
Half-facets on the boundary of the domain have no siblings. The half-facets
are \emph{half-edges} and \emph{half-faces} in 2-D and 3-D, respectively.
We identify each half-facet by a two tuple: the element ID and a local
facet ID within the element. In 2-D, we store the element connectivity,
sibling half-edges, and a mapping from each node to an incident half-edge.
In 3-D, we store the element connectivity, sibling half-faces, and
a mapping from each node to an incident half-face. This data structure
allows us to do neighborhood queries for a node in constant time (provided
the valance is bounded). For additional information about the data
structure, see \cite{DREJTAHF2014}.

\subsubsection{Solution of Weighted Least Squares Problems}

The technique for solving this least squares problem (\ref{eq:average_eq})
is QR factorization with column pivoting. To further improve the condition
number of $\vec{A}$, we follow the idea of \cite{JiaoWang10RHO,jiao2008consistent}
by scaling the columns of the matrix $\vec{W}\vec{A}$ and solving
the following problem instead:
\begin{equation}
\underset{d}{\min}\left\Vert \vec{W}\vec{A}\vec{S}\vec{d}-\vec{W}\overline{\vec{u}}\right\Vert _{2}.\label{eq:scaled_WLS}
\end{equation}
Here $\vec{W}$ is the weighting matrix in (\ref{eq:WLS}), $\vec{d}\equiv\vec{S}^{-1}\vec{V}$
and $\vec{S}=\mathrm{diag}\left(1/\left\Vert \tilde{\vec{a}}_{1}\right\Vert _{2},1/\left\Vert \tilde{\vec{a}}_{2}\right\Vert _{2},\dots,1/\left\Vert \tilde{\vec{a}}_{n}\right\Vert _{2}\right)$,
where $\tilde{\vec{a}}_{i}$ is the $i$th column vector of $\vec{W}\vec{A}$.
We perform reduced QR factorization with column pivoting to the matrix
$\vec{W}\vec{A}\vec{S}$:
\begin{equation}
\vec{W}\vec{A}\vec{S}\vec{E}=\vec{Q}\vec{R}.
\end{equation}
Here $\vec{E}$ is chosen so that the diagonal of $\vec{R}$ is in
decreasing order. If $\vec{W}\vec{A}\vec{S}$ has full rank, then
its pseudoinverse is
\begin{equation}
\left(\vec{W}\vec{A}\vec{S}\right)^{+}=\vec{E}\vec{R}^{-1}\vec{Q}^{T}.
\end{equation}
Otherwise, the pseudoinverse is computed as
\begin{equation}
\left(\vec{W}\vec{A}\vec{S}\right)^{+}=\vec{E}_{1:k,1:r}\vec{R}_{1:r,1:r}^{-1}(\vec{Q}_{1:m,1:r})^{T},
\end{equation}
where $r$ the numerical rank of $\vec{R}$. In this way, we can truncate
the higher order terms in $\vec{W}\vec{A}\vec{S}$ for best-possible
accuracy whenever possible.

\section{\label{sec:Analysis-of-LS-WENO}Accuracy and Stability of WLS-ENO}

In this section, we analyze the accuracy of WLS-ENO, as well as its
stability in the context of solving hyperbolic conservation laws.

\subsection{\label{sub:Accuracy}Accuracy}

First, we analyze the WLS-ENO schemes and show that they can achieve
the expected order of accuracy for smooth functions.

\label{thm:accuracy}Given a mesh with a smooth function $f$. Let
$\vec{W}$ be a diagonal matrix containing all the weights for the
cells. $\vec{A}$ and $\vec{S}$ are the matrices in (\ref{eq:scaled_WLS}).
Suppose the cell average of $f$ is approximated with an error $\mathcal{O}(h^{p})$
and the matrix $\vec{W}\vec{A}\vec{S}$ has a bounded condition number.
The degree-($p-1$) cell average weighed least squares fitting approximates
$q$th order derivatives of function $f$ to $\mathcal{O}(h^{p-q})$.

For simplicity, we only prove the theorem in 2-D. The analysis also
applies to 1-D and 3-D.

The 2-D Taylor series expansion about the point $(x_{i},y_{i})$ reads
\begin{equation}
f(x,y)=\sum_{q=0}^{p-1}\sum_{j,k\geq0}^{j+k=q}f_{jk}\frac{(x-x_{i})^{j}(y-y_{i})^{k}}{j!k!}+\mathcal{O}\left(\left\Vert \delta\vec{x}\right\Vert ^{p}\right),
\end{equation}
where $\delta\vec{x}=\left[x-x_{i},y-y_{i}\right]^{T}$. The cell
average of $f(x,y)$ over some cell $\tau_{i}$ can be written as
\begin{equation}
\frac{1}{\left|\tau_{i}\right|}\iint_{\tau_{i}}f(x,y)\,dxdy=\frac{1}{\left|\tau_{i}\right|}\sum_{q=0}^{p-1}\sum_{j,k\geq0}^{j+k=q}\frac{f_{jk}}{j!k!}\iint_{\tau_{i}}(x-x_{i})^{j}(y-y_{i})^{k}\,dxdy+\mathcal{O}\left(\left\Vert \delta\vec{x}\right\Vert ^{p}\right).
\end{equation}
Let $\vec{v}$ denotes the exact derivatives of function $f$, $\tilde{\vec{v}}$
the numerical solution from the WLS fitting. Let $\vec{r}=\overline{\vec{u}}-\vec{A}\vec{v}$.
By assumption, each component of $\vec{r}$ is $\mathcal{O}\left(\left\Vert \delta\vec{x}\right\Vert ^{p}\right)$.
The error of coefficients has the relationship $\vec{W}\vec{A}\left(\tilde{\vec{v}}-\vec{v}\right)\approx\vec{W}\vec{r}$.
The error of $\vec{d}$ can then be written as
\begin{equation}
\vec{W}\vec{A}\vec{S}\delta\vec{d}\approx\vec{W}\vec{r}.
\end{equation}
By solving this least squares problem, we have $\delta\vec{d}=\left(\vec{W}\vec{A}\vec{S}\right)^{+}\vec{W}\vec{r}$.
Since the function $f$ is smooth, all the diagonal entries in $\vec{W}$
are $\mathcal{O}(1/\left\Vert \delta\vec{x}\right\Vert ^{2})$. Under
the assumption that $\vec{W}\vec{A}\vec{S}$ has a bounded condition
number $\kappa$, all the component of $\delta\vec{d}$ are $\mathcal{O}(\kappa\left\Vert \delta\vec{x}\right\Vert ^{p-2})$.
For a $q$th order partial derivative of function $f$, the corresponding
column in $\vec{W}\vec{A}$ is $\mathcal{O}(\left\Vert \delta\vec{x}\right\Vert ^{q-2})$,
so is the 2-norm of the column. Therefore, the $q$th order derivatives
of function $f$ are approximated to $\mathcal{O}(h^{p-q})$.

From this theorem, we can conclude that the degree-$(p-1)$ WLS-ENO
reconstruction delivers a $p$th order accurate reconstruction for
smooth functions.

\subsection{\label{sub:Stability-of-WLS-ENO}Stability of WLS-ENO for Hyperbolic
Conservation Laws}

For the WLS-ENO to be practically useful, it is important that it
enables a stable discretization for hyperbolic conservation laws,
when coupled with a proper time-integration scheme. In the following,
we analyze two fifth-order WLS-ENO schemes for a model problem in
1-D, based on a modified von Neumann stability analysis. We show that
the WLS-ENO scheme has a larger stability region than the WENO scheme
on structured meshes.

\subsubsection{Model Problem in 1-D}

We consider one dimensional wave equation
\begin{equation}
u_{t}+u_{x}=0,\hspace{1em}x\in[0,1],\hspace{1em}t>0
\end{equation}
with the periodic boundary condition
\begin{equation}
u(x,0)=u^{0}(x),\hspace{1em}x\in[0,1].
\end{equation}
Suppose we have a structured grid $0=x_{\frac{1}{2}}<x_{\frac{3}{2}}<\cdots<x_{N+\frac{1}{2}}=1$
with $x_{i}=i\Delta x$ and $\Delta x=1/N$. For the finite volume
method, we integrate the above wave equation and divide it by the
length of the cell, and obtain
\begin{equation}
\frac{d\overline{u}(x_{i},t)}{dt}=-\frac{1}{\Delta x_{i}}\left(f\left(u\left(x_{i+\frac{1}{2}},t\right)\right)-f\left(u\left(x_{i-\frac{1}{2}},t\right)\right)\right),\label{eq:1D FV}
\end{equation}
where
\begin{equation}
\overline{u}\left(x_{i},t\right)=\frac{1}{\Delta x_{i}}\int_{x_{i-\frac{1}{2}}}^{x_{i+\frac{1}{2}}}u(x,t)\,dx.
\end{equation}
We approximate equation (\ref{eq:1D FV}) by the following conservative
scheme
\begin{equation}
\frac{d\overline{u}(x_{i},t)}{dt}=-\frac{1}{\Delta x_{i}}\left(\hat{f}_{i+\frac{1}{2}}-\hat{f}_{i-\frac{1}{2}}\right),
\end{equation}
where the numerical flux $\hat{f}_{i+\frac{1}{2}}$ is replaced by
the Lax-Friedrichs flux (\ref{eq:flux}), with $\alpha=\mathrm{max}_{u}\left|f'(u)\right|$.
For this particular problem, $\hat{f}_{i+\frac{1}{2}}-\hat{f}_{i-\frac{1}{2}}=u_{i+\frac{1}{2}}^{-}-u_{i-\frac{1}{2}}^{-}$.

For time integration, we use a TVD Runge-Kutta method \cite{gottlieb1998total}.
An $n$-stage Runge-Kutta method for the ODE $u_{t}=L(u)$ has the
general form of
\begin{align}
k_{0} & =u(t),\\
k_{i} & =\sum_{j=0}^{i-1}\left(\alpha_{ij}k_{j}+\beta_{ij}\Delta tL(k_{j})\right),\qquad i=1,\dots,n,\label{eq:runge-kutta}
\end{align}
where $k_{i}$ denotes the intermediate solution after the $i$th
stage, and $u(t+\Delta t)=k_{n}$. A Runge-Kutta method is total variation
diminishing (TVD) if all the coefficients $\alpha_{ij}$ and $\beta_{ij}$
are nonnegative. The CFL coefficient of such a scheme is given by
\begin{equation}
c=\underset{i,k}{\min}\{\alpha_{ik}/\beta_{ik}\}.
\end{equation}
Specifically, we use the third-order TVD Runge-Kutta scheme, given
by
\begin{align}
k_{1} & =u+\Delta tL\left(u\right),\\
k_{2} & =\frac{3}{4}u+\frac{1}{4}k_{1}+\frac{1}{4}\Delta tL\left(k_{1}\right),\\
k_{3} & =\frac{1}{3}u+\frac{2}{3}k_{2}+\frac{2}{3}\Delta tL(k_{2}),\label{eq:3rd-order-Runge-Kutta}
\end{align}
for which the CFL coefficient is $c=1$.

\subsubsection{von Neumann Stability Analysis}

Based on the von Neumann stability analysis, the semi-discrete solution
can be expressed by a discrete Fourier series
\begin{equation}
\overline{u}_{j}(t)=\sum_{k=-N/2}^{N/2}\hat{u}_{k}(t)e^{i\omega_{k}j\Delta x},\hspace{1em}\omega_{k}\in R.
\end{equation}
By the superposition principle, we can only use one term in the series
for analysis
\begin{equation}
\overline{u}_{j}(t)=\hat{u}_{k}(t)e^{ij\theta_{k}},\hspace{1em}\theta_{k}=\omega_{k}\Delta x,\label{eq:von Neumann}
\end{equation}
where $k=-N/2,\dots,N/2$. We assume that the numerical flux can be
written in the following form
\begin{equation}
\hat{f}_{i+\frac{1}{2}}-\hat{f}_{i-\frac{1}{2}}=u_{i+\frac{1}{2}}^{-}-u_{i-\frac{1}{2}}^{-}=z\left(\theta_{k}\right)\overline{u}_{i},
\end{equation}
where the complex function $z\left(\theta_{k}\right)$ is the Fourier
symbol. 

Let $\overline{u}_{i}^{n}=\overline{u}_{i}(t^{n})$ be the numerical
solution at time level $t^{n}=n\Delta t$. We define the amplification
factor $g$ by inserting (\ref{eq:von Neumann}) into the fully-discrete
system and obtain
\begin{equation}
\overline{u}_{i}^{n+1}=g(\hat{z}_{k})\overline{u}_{i}^{n},\hspace{1em}\hat{z}_{k}=-\sigma z\left(\theta_{k}\right),\hspace{1em}k=-N/2\dots N/2,
\end{equation}
where $\sigma=\Delta t/\Delta x$. Therefore, the linear stability
domain of an explicit time-stepping scheme is $S_{t}=\left\{ \hat{z}:\left|g\left(\hat{z}\right)\right|\leq1\right\} $.
Also, we define the discrete spectrum $S$ of a spatial discretization
scheme
\begin{equation}
S=\left\{ -z\left(\theta_{k}\right):\theta_{k}\in0,\Delta\theta,2\Delta\theta,\dots,2\pi\right\} ,\hspace{1em}\Delta\theta=2\pi\Delta x.
\end{equation}
The stability limit is thus the largest CFL number $\widetilde{\sigma}$
such that the rescaled spectrum $\widetilde{\sigma}S$ lies inside
the stability domain
\begin{equation}
\widetilde{\sigma}S\in S_{t}.
\end{equation}
For the third-order Runge-Kutta scheme, the amplification factor is
given by
\begin{equation}
g\left(\widetilde{z}\right)=1+\widetilde{z}+\frac{1}{2}\widetilde{z}^{2}+\frac{1}{6}\widetilde{z}^{3}.
\end{equation}
To determine the boundary of the stability domain $\partial S_{t}=\left\{ \widetilde{z}:\left|g(\widetilde{z})\right|=1\right\} $,
we set $g(\widetilde{z})=e^{i\phi}$ and solve the following equation
\begin{equation}
\widetilde{z}^{3}+3\widetilde{z}^{2}+6\widetilde{z}+6-6e^{i\phi}=0.
\end{equation}
The stability region is given by the solid blue curves in Figure~\ref{fig:The-rescaled-spectrum}.

\subsubsection{Fifth-Order WLS-ENO Scheme with Five Cells}

Let us first consider a fifth-order scheme for a cell $\tau_{i}$,
with five cells
\begin{equation}
S_{i}=\left\{ \tau_{i-2},\tau_{i-1},\tau_{i},\tau_{i+1},\tau_{i+2}\right\} ,
\end{equation}
where $\tau_{i}=\left[x_{i-\frac{1}{2}},x_{i+\frac{1}{2}}\right]$.
This results in the linear system (\ref{eq:average_eq})
\begin{equation}
\vec{A}=\left(\begin{array}{ccccc}
1 & -\frac{5\Delta x}{2} & \frac{19\Delta x^{2}}{6} & -\frac{65\Delta x^{3}}{24} & \frac{211\Delta x^{4}}{120}\\
1 & -\frac{3\Delta x}{2} & \frac{7\triangle x^{2}}{6} & -\frac{15\Delta x^{3}}{24} & \frac{31\Delta x^{4}}{120}\\
1 & -\frac{\Delta x}{2} & \frac{\Delta x^{2}}{6} & -\frac{\Delta x^{3}}{24} & \frac{\Delta x^{4}}{120}\\
1 & \frac{\Delta x}{2} & \frac{\Delta x^{2}}{6} & \frac{\Delta x^{3}}{24} & \frac{\Delta x^{4}}{120}\\
1 & \frac{3\Delta x}{2} & \frac{7\Delta x^{2}}{6} & \frac{15\Delta x^{3}}{24} & \frac{31\Delta x^{4}}{120}
\end{array}\right).\label{eq:5 cell 5th order matrix}
\end{equation}
Since $\vec{A}$ is nonsingular, the weights do not affect the solution.
The solution is given by 
\begin{equation}
u_{i+\frac{1}{2}}^{-}=\frac{2}{60}\overline{u}_{i-2}-\frac{13}{60}\overline{u}_{i-1}+\frac{47}{60}\overline{u}_{i}+\frac{27}{60}\overline{u}_{i+1}-\frac{3}{60}\overline{u}_{i+2}.
\end{equation}
This is equivalent to the fifth-order WENO scheme without nonlinear
weights. Also, the flux reads
\begin{equation}
u_{i+\frac{1}{2}}^{-}-u_{i-\frac{1}{2}}^{-}=-\frac{2}{60}\overline{u}_{i-3}+\frac{15}{60}\overline{u}_{i-2}-\frac{60}{60}\overline{u}_{i-1}+\frac{20}{60}\overline{u}_{i}+\frac{30}{60}\overline{u}_{i+1}-\frac{3}{60}\overline{u}_{i+2}.
\end{equation}
Substituting it into (\ref{eq:von Neumann}), we get
\begin{equation}
z\left(\theta_{k}\right)=\frac{16}{15}\sin^{6}\left(\frac{\theta_{k}}{2}\right)+i\left(-\frac{1}{6}\sin\left(2\theta_{k}\right)+\frac{4}{3}\sin\left(\theta_{k}\right)+\frac{16}{15}\sin^{5}\left(\frac{\theta_{k}}{2}\right)\cos\left(\frac{\theta_{k}}{2}\right)\right).
\end{equation}
Because the eigenvalues $-z\left(\theta_{k}\right)$ on the discrete
spectrum go in clockwise order as $\theta_{k}$ increase from 0 to
$2\pi$, we reverse its order by replacing $\theta_{k}$ by $2\pi-\theta_{k}$,
so that it is consistent with the standard polar form in the complex
plane
\begin{equation}
z\left(\theta_{k}\right)=\frac{16}{15}\sin^{6}\left(\frac{\theta_{k}}{2}\right)+i\left(\frac{1}{6}\sin\left(2\theta_{k}\right)-\frac{4}{3}\sin\left(\theta_{k}\right)-\frac{16}{15}\sin^{5}\left(\frac{\theta_{k}}{2}\right)\cos\left(\frac{\theta_{k}}{2}\right)\right).
\end{equation}
The discrete spectrum is shown in Figure~\ref{fig:The-rescaled-spectrum}(left).

Given the spectrum and the stability domain, the CFL number of this
scheme can be computed by finding the largest rescaling parameter
$\sigma$, so that the rescaled spectrum still lies in the stability
domain. Using interval bisection, we find that the CFL number is $\sigma=1.44$
for the fifth-order five-cell scheme.

\subsubsection{Fifth-Order WLS-ENO Scheme with Seven Cells}

\begin{figure}
\begin{centering}
\begin{minipage}[t]{0.49\columnwidth}%
\begin{center}
\includegraphics[width=1\textwidth]{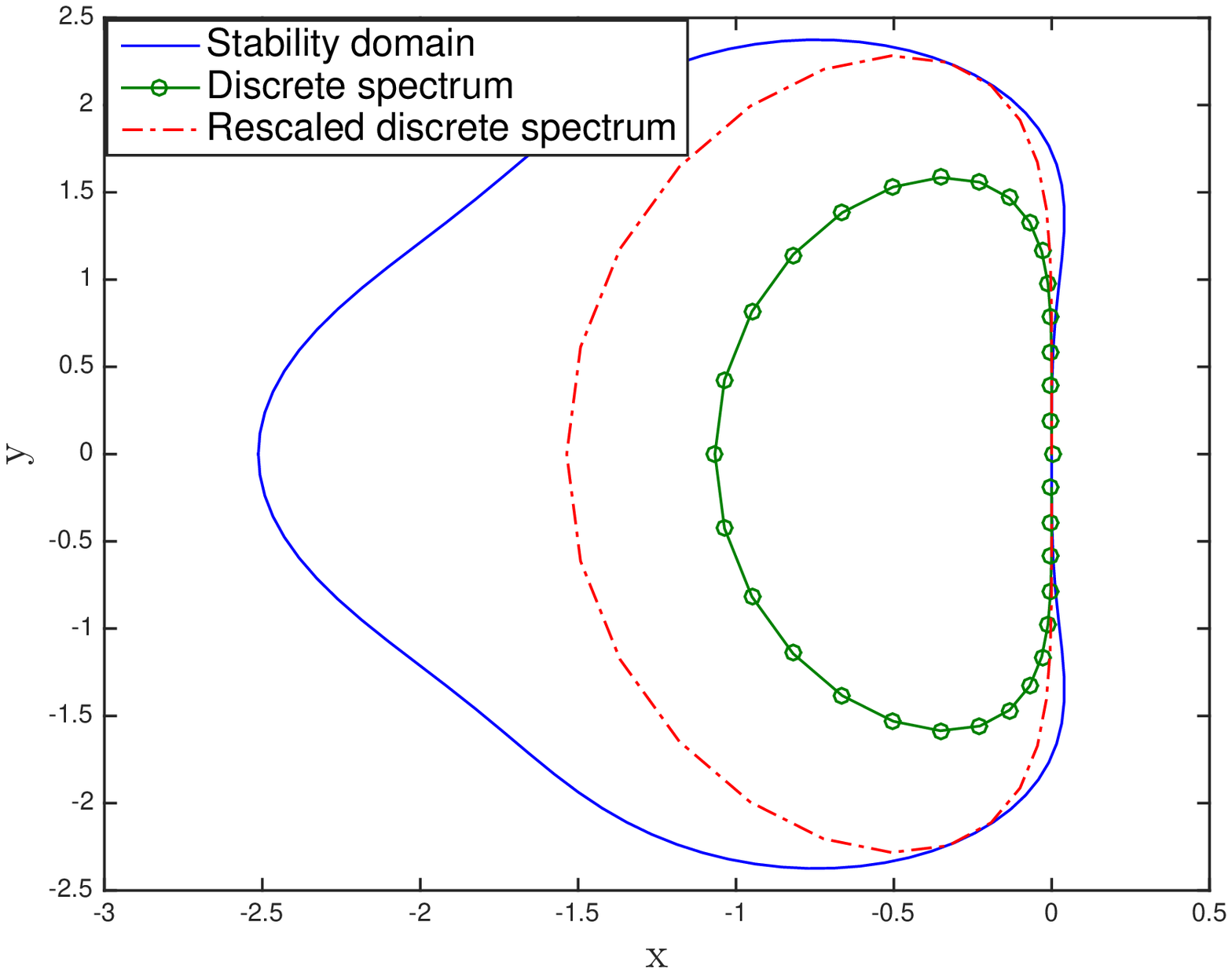}
\par\end{center}%
\end{minipage}\hfill%
\begin{minipage}[t]{0.49\columnwidth}%
\begin{center}
\includegraphics[width=1\textwidth]{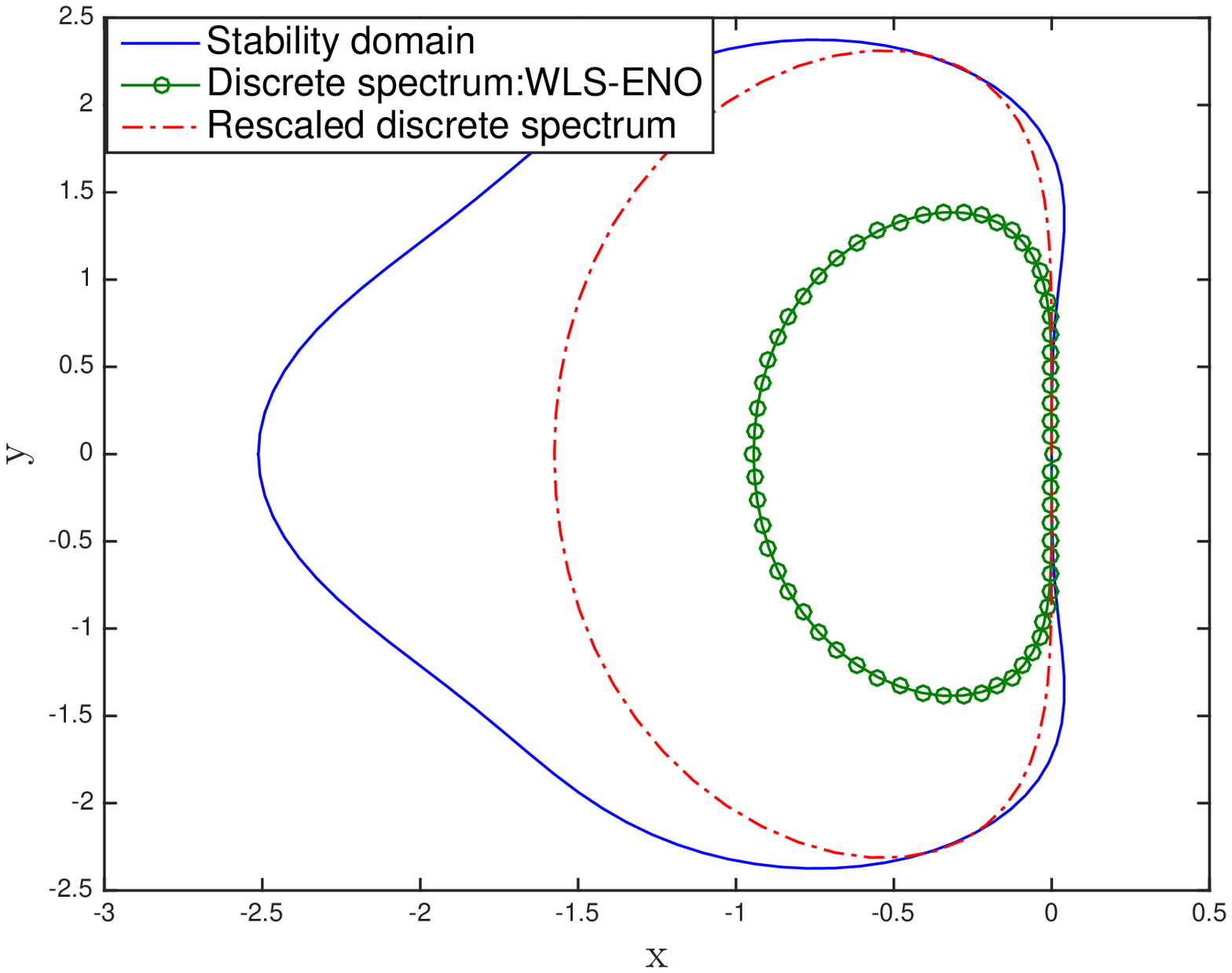}
\par\end{center}%
\end{minipage}
\par\end{centering}

\protect\caption{\label{fig:The-rescaled-spectrum}Rescaled spectrum and stability
domains of fifth order WLS-ENO with five (left) and seven (right)
cells.}
\end{figure}

Next, let us consider a least-squares fitting for $\tau_{i}$ using
seven cells
\begin{equation}
S_{i}=\left\{ \tau_{i-3},\tau_{i-2},\tau_{i-1},\tau_{i},\tau_{i+1},\tau_{i+2},\tau_{i+3}\right\} .
\end{equation}
The linear system is given by 
\begin{equation}
\vec{A}=\left(\begin{array}{ccccc}
1 & -\frac{7\Delta x}{2} & \frac{37\Delta x^{2}}{6} & -\frac{175\Delta x^{3}}{24} & \frac{781\Delta x^{4}}{120}\\
1 & -\frac{5\Delta x}{2} & \frac{19\Delta x^{2}}{6} & -\frac{65\Delta x^{3}}{24} & \frac{211\Delta x^{4}}{120}\\
1 & -\frac{3\Delta x}{2} & \frac{7\triangle x^{2}}{6} & -\frac{15\Delta x^{3}}{24} & \frac{31\Delta x^{4}}{120}\\
1 & -\frac{\Delta x}{2} & \frac{\Delta x^{2}}{6} & -\frac{\Delta x^{3}}{24} & \frac{\Delta x^{4}}{120}\\
1 & \frac{\Delta x}{2} & \frac{\Delta x^{2}}{6} & \frac{\Delta x^{3}}{24} & \frac{\Delta x^{4}}{120}\\
1 & \frac{3\Delta x}{2} & \frac{7\Delta x^{2}}{6} & \frac{15\Delta x^{3}}{24} & \frac{31\Delta x^{4}}{120}\\
1 & \frac{5\Delta x}{2} & \frac{19\Delta x^{2}}{6} & \frac{65\Delta x^{3}}{24} & \frac{211\Delta x^{4}}{120}
\end{array}\right).\label{eq:7 cell 5th order matrix}
\end{equation}
If $\epsilon=0$, from (\ref{eq:smoothness indicator}) and (\ref{eq:weight}),
we obtain the weights
\begin{equation}
w_{1}\approx\frac{1}{9\Delta x^{2}},\hspace{1em}w_{2}\approx\frac{1}{4\Delta x^{2}},\hspace{1em}w_{3}\approx\frac{1}{\Delta x^{2}},\hspace{1em}w_{4}\approx\frac{1.5}{\Delta x^{2}},\hspace{1em}w_{5}\approx\frac{1}{\Delta x^{2}},\hspace{1em}w_{6}\approx\frac{1}{4\Delta x^{2}},\hspace{1em}w_{7}\approx\frac{1}{9\Delta x^{2}}.
\end{equation}
Solving the weighted least squares system, we obtain the following
scheme
\begin{align*}
u_{i+\frac{1}{2}}^{-}= & \frac{1226983}{9489680}\overline{u}_{i-3}-\frac{963431}{47447340}\overline{u}_{i-4}-\frac{13515169}{94894680}\overline{u}_{i-1}+\frac{66771}{87380}\overline{u}_{i}\\
 & +\frac{38388551}{94894680}\overline{u}_{i+1}-\frac{93404}{11861835}\overline{u}_{i+2}-\frac{348299}{31631560}\overline{u}_{i+3}.
\end{align*}
Therefore, the flux reads
\begin{align*}
u_{i+\frac{1}{2}}^{-}-u_{i-\frac{1}{2}}^{-}= & -\frac{1226983}{94894680}\overline{u}_{i-4}+\frac{630769}{18978936}\overline{u}_{i-3}+\frac{3862769}{31631560}\overline{u}_{i-2}-\frac{17205695}{18978936}\overline{u}_{i-1}\\
 & +\frac{6824951}{18978936}\overline{u}_{i}+\frac{13045261}{31631560}\overline{u}_{i+1}+\frac{59533}{18978936}\overline{u}_{i+2}-\frac{348299}{31631560}\overline{u}_{i+3.}
\end{align*}
Substituting it into (\ref{eq:von Neumann}), we get
\begin{equation}
z\left(\theta_{k}\right)=p\left(\theta_{k}\right)+iq\left(\theta_{k}\right),
\end{equation}
where
\begin{equation}
p\left(\theta_{k}\right)=-0.0129\cos\left(4\theta_{k}\right)+0.0222\cos\left(3\theta_{k}\right)+0.1253\cos\left(2\theta_{k}\right)-0.4942\cos\left(\theta_{k}\right)+0.3596,
\end{equation}
and
\begin{equation}
q\left(\theta_{k}\right)=0.0129\sin\left(4\theta_{k}\right)-0.0442\sin\left(3\theta_{k}\right)-0.1190\sin\left(2\theta_{k}\right)+1.3190\sin\left(\theta_{k}\right).
\end{equation}
The discrete spectrum is shown in Figure~\ref{fig:The-rescaled-spectrum}(right).
The CFL number is computed in the same way as the above scheme. For
this scheme with seven cells, we obtain $\sigma=1.67$. Therefore,
the seven-cell fifth-order WLS-ENO scheme has a larger stability region
than the five-cell counterpart, or equivalently the fifth-order WENO
scheme.

\section{\label{sec:Numerical-Results}Numerical Results}

In this section, we present some numerical experiments of WLS-ENO
over both structured and unstructured meshes in 1-D, 2-D, and 3-D,
and compare it against WENO schemes when applicable. For all the PDEs,
we will use third-order TVD Runge-Kutta for time integration.

\subsection{1-D Results}

We first show some results in 1-D, for the reconstruction of a piecewise
smooth function as well as the solutions of PDEs, including a linear
wave equation, Burgers' equation, and Euler equations.

\subsubsection{Reconstruction of Discontinuous Functions}

\begin{figure}
\begin{minipage}[t]{0.49\columnwidth}%
\begin{center}
\includegraphics[width=1\textwidth]{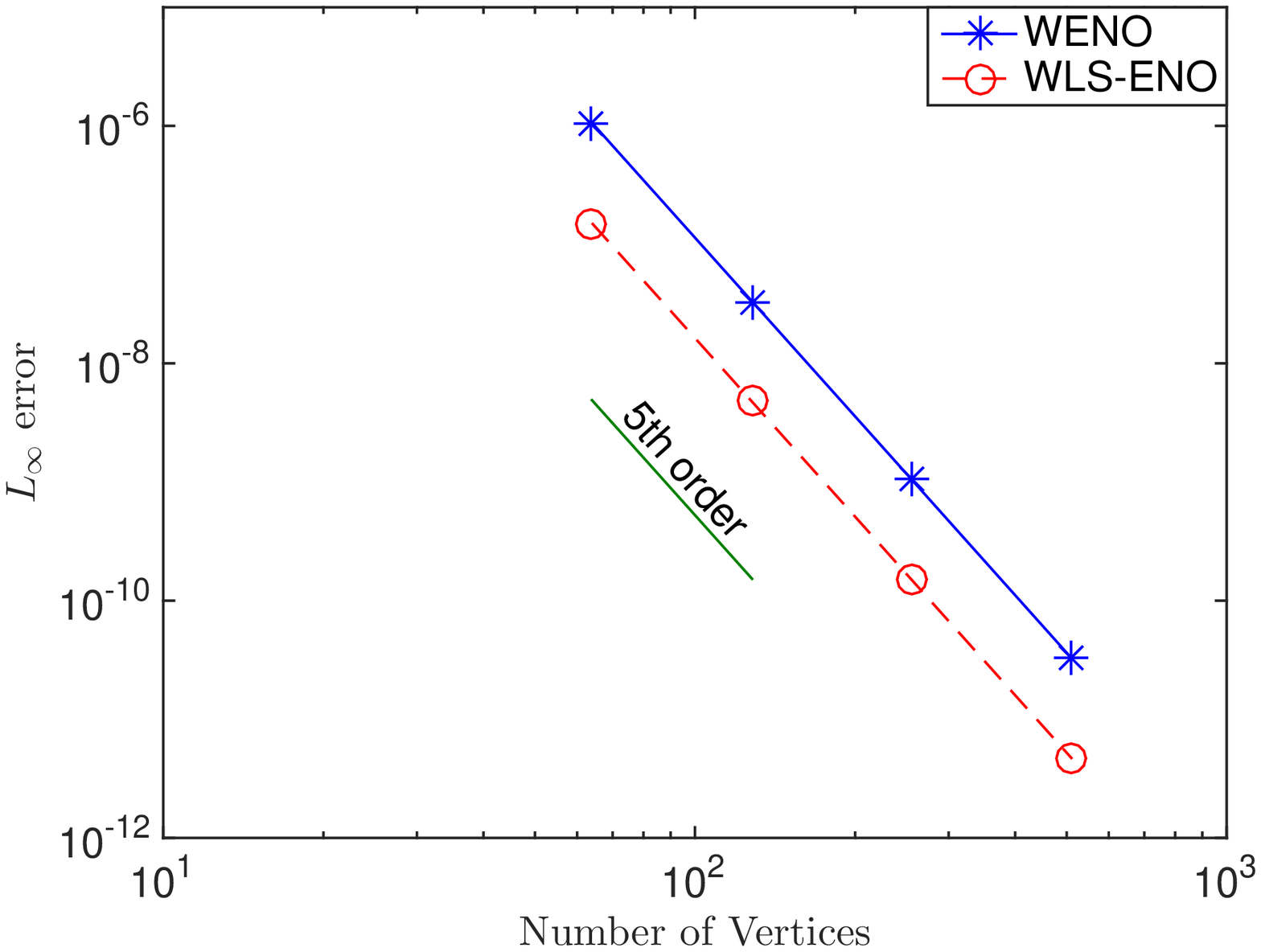}
\par\end{center}

\protect\caption{\label{fig:1d-convergence-smooth}Convergence of fifth order WENO
and WLS-ENO away from discontinuity.}
\end{minipage}\hfill%
\begin{minipage}[t]{0.49\columnwidth}%
\begin{center}
\includegraphics[width=1\textwidth]{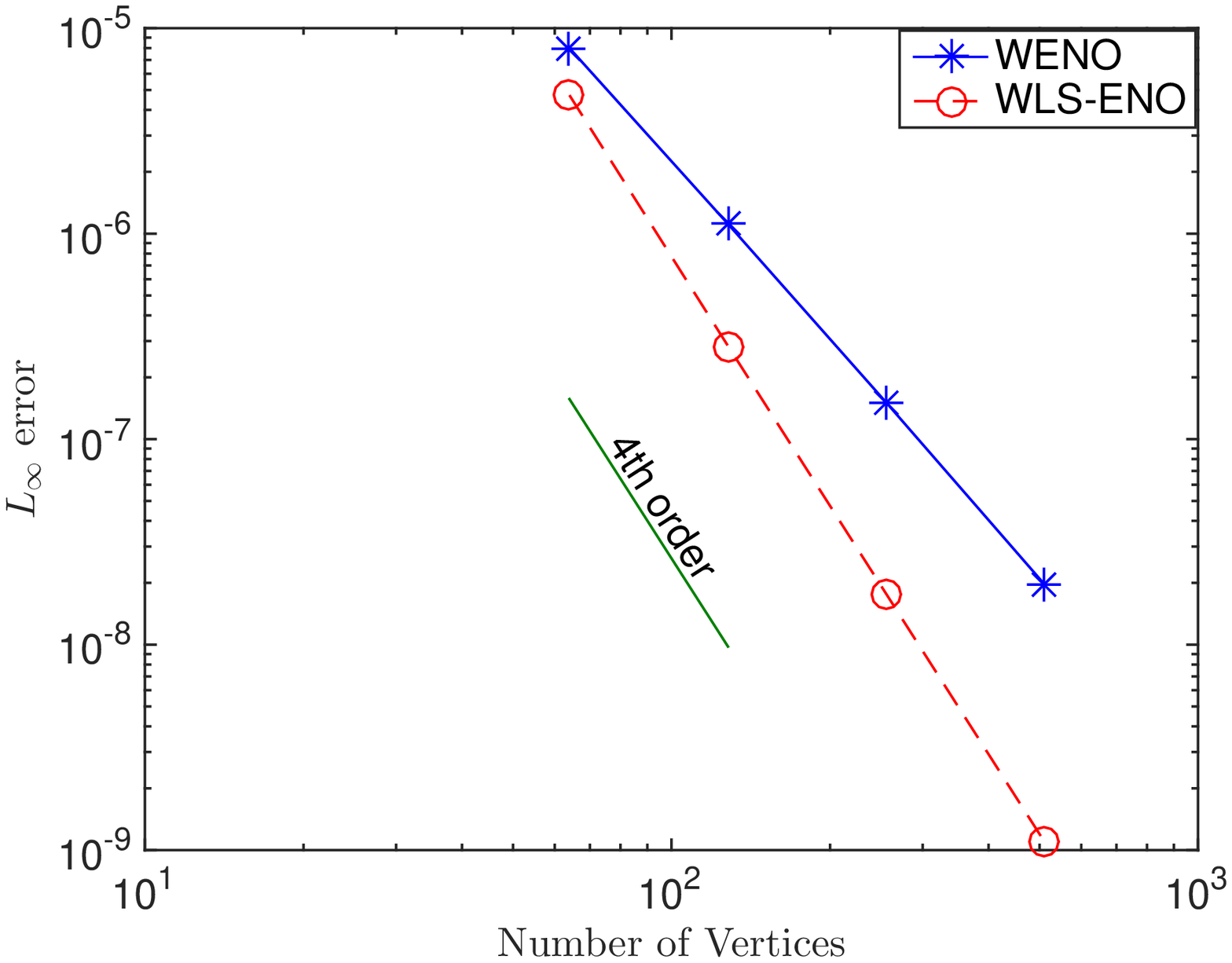}
\par\end{center}

\begin{center}
\protect\caption{\label{fig:1d-conv-singularity}Convergence of fifth order WENO and
WLS-ENO near discontinuity.}

\par\end{center}%
\end{minipage}
\end{figure}

We first test WLS-ENO for the reconstruction of a 1-D discontinuous
function, given by 
\begin{equation}
v(x)=\begin{cases}
\begin{array}{c}
\sin(\pi x)\\
\cos(\pi x)
\end{array} & \begin{array}{c}
0\leq x\leq0.6\\
0.6<x\leq1
\end{array}\end{cases}.
\end{equation}
This function is discontinuous at $x=0.6$ but smooth everywhere else
within the interval $[0,1]$. We performed grid convergence study
under grid refinement, starting from an equidistant grid with 32 grid
cells. For the WLS-ENO, we used degree-four polynomials over seven-cell
stencil, which based on our theory should deliver fifth-order accuracy
in smooth regions and fourth-order accuracy near discontinuities.
As a point of reference, we also perform the reconstruction using
the fifth-order WENO scheme, which is fifth-order accurate in smooth
regions and third-order accurate near discontinuities. Figure~\ref{fig:1d-convergence-smooth}
shows the $L_{\infty}$-norm error for the reconstructed values at
the grid points that are one cell away from the discontinuity. It
can be seen that both WLS-ENO and WENO delivered fifth order accuracy,
but WLS-ENO is more accurate. When including the grid points near
discontinuities, as can be seen in Figure~\ref{fig:1d-conv-singularity},
WLS-ENO achieved fourth order convergence, whereas the fifth-order
WENO reduced to third order, as predicted by their respective theoretical
analyses.

\subsubsection{1-D Wave Equation}

\begin{figure}
\begin{centering}
\begin{minipage}[t]{0.49\columnwidth}%
\begin{center}
\includegraphics[width=1\textwidth]{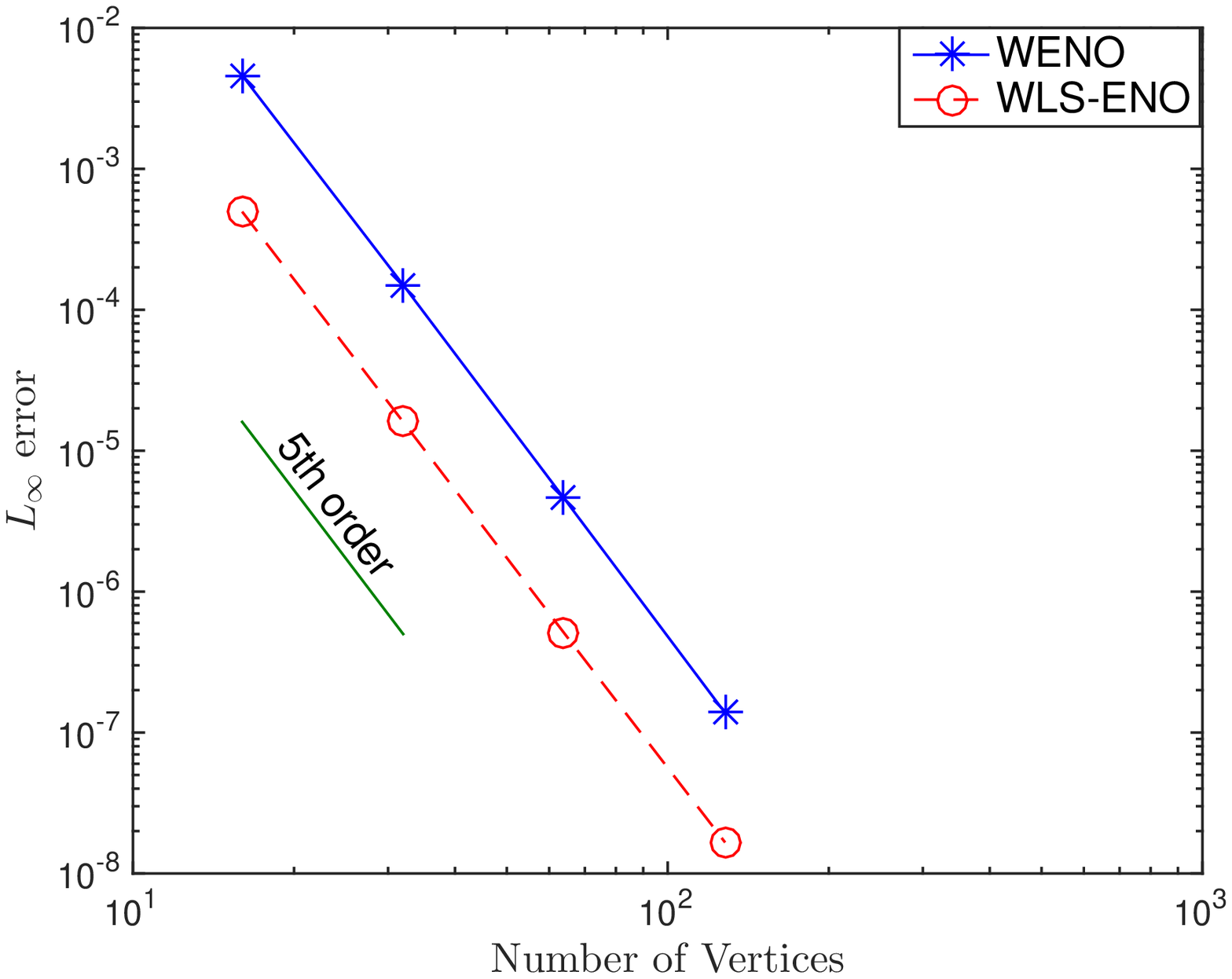}\\
{\small{}(a) Errors on uniform grids.}
\par\end{center}{\small \par}%
\end{minipage}\hfill %
\begin{minipage}[t]{0.49\columnwidth}%
\begin{center}
\includegraphics[width=1\textwidth]{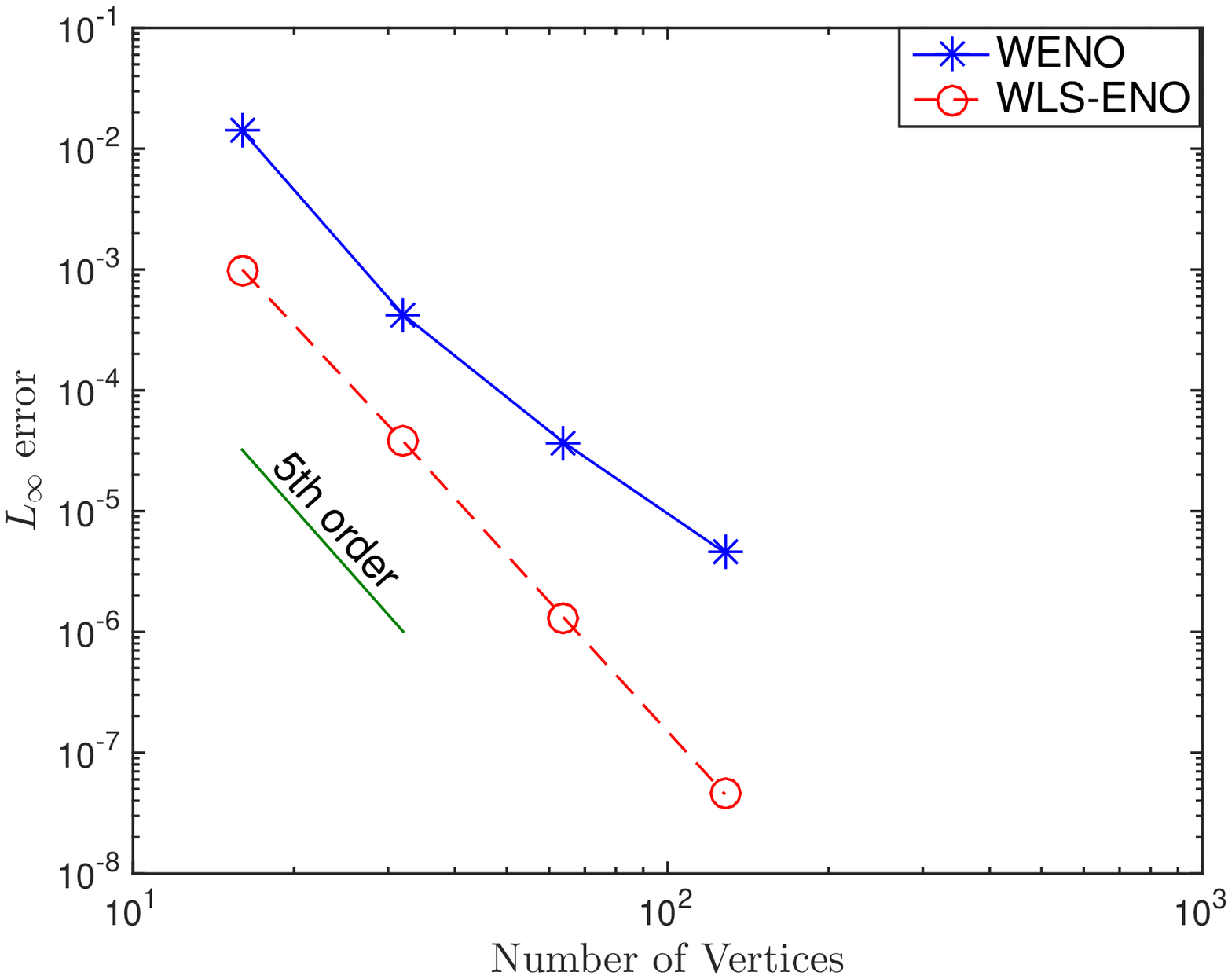}\\
{\small{}(b) Errors on non-uniform grids.}
\par\end{center}{\small \par}%
\end{minipage}
\par\end{centering}

\centering{}\protect\caption{\label{fig:Convergence-linwaveeq}Convergence of fifth order WENO
and WLS-ENO for linear wave equation at $t=1$ on (a) uniform and
(b) non-uniform grids.}
\end{figure}

To test the effectiveness of WLS-ENO for solving hyperbolic PDEs,
we first solve a simple linear wave equation
\begin{equation}
u_{t}+u_{x}=0,\qquad-1\leq x\leq1,
\end{equation}
with periodic boundary conditions. Similar to the reconstruction problem,
we use WLS-ENO with degree-four polynomials over seven cells. To assess
the accuracy for smooth solutions, consider the smooth initial condition
\begin{equation}
u(x,0)=\sin(\pi x),
\end{equation}
for which the solution remains smooth over time. We assess the order
of accuracy of the solutions at $t=1$ under grid refinement, and
compare the errors against the fifth-order WENO. Figure~\ref{fig:Convergence-linwaveeq}
shows the results for uniform and non-uniform grids. For uniform grids,
similar to the results of reconstruction, both WLS-ENO and WENO delivered
fifth order convergence under grid refinement and solution of WLS-ENO
scheme is more accurate. For non-uniform grids, we used WENO reconstruction
described in \cite{ShuWeno98}, which converged at a slower rates,
and was about one order of magnitude less accurate than WLS-ENO on
the finest grid. 

\begin{figure}
\begin{centering}
\begin{minipage}[t]{0.49\columnwidth}%
\begin{center}
\includegraphics[width=1\textwidth]{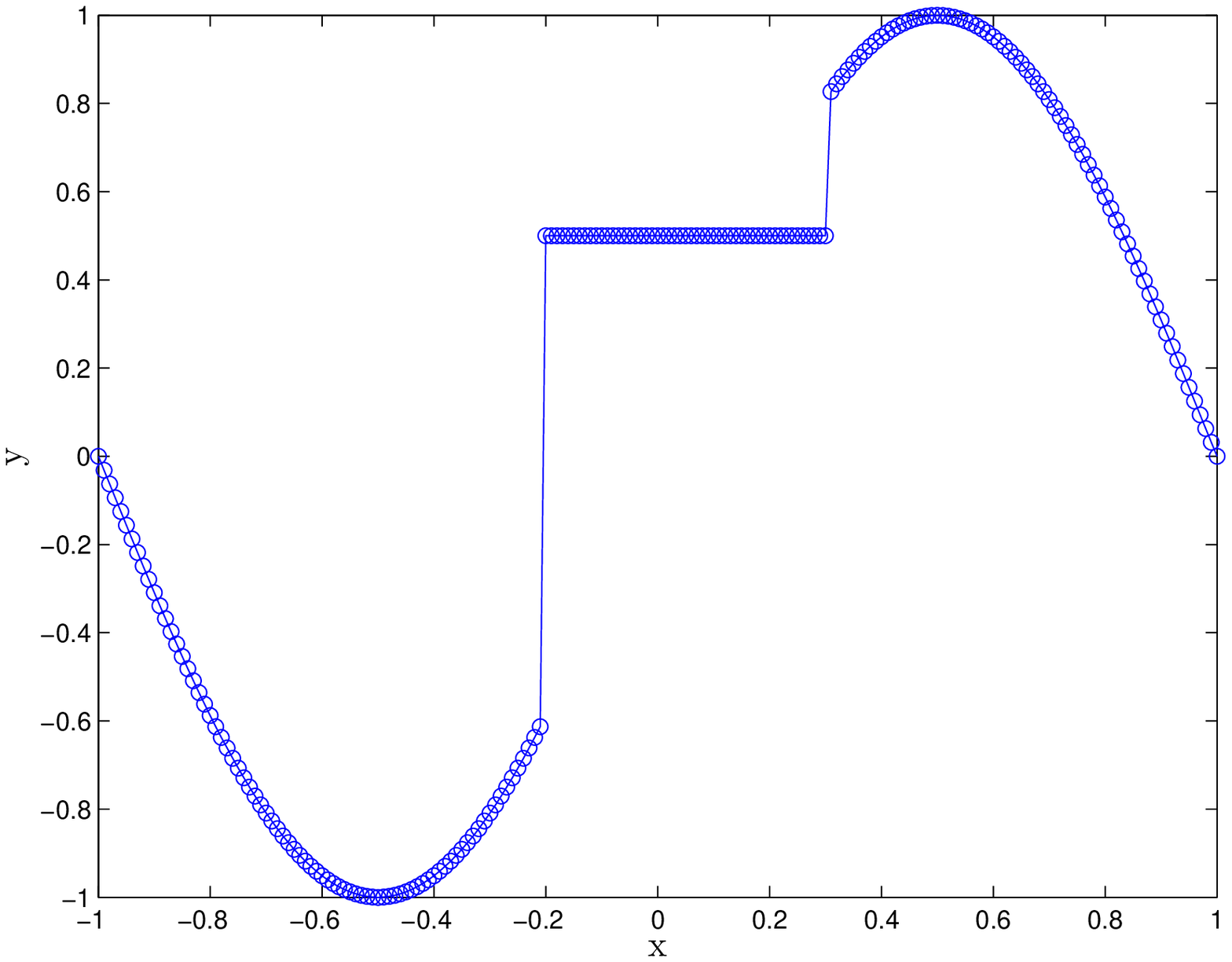}\\
{\small{}(a) Initial condition.}
\par\end{center}{\small \par}%
\end{minipage}\hfill %
\begin{minipage}[t]{0.49\columnwidth}%
\begin{center}
\includegraphics[width=1\textwidth]{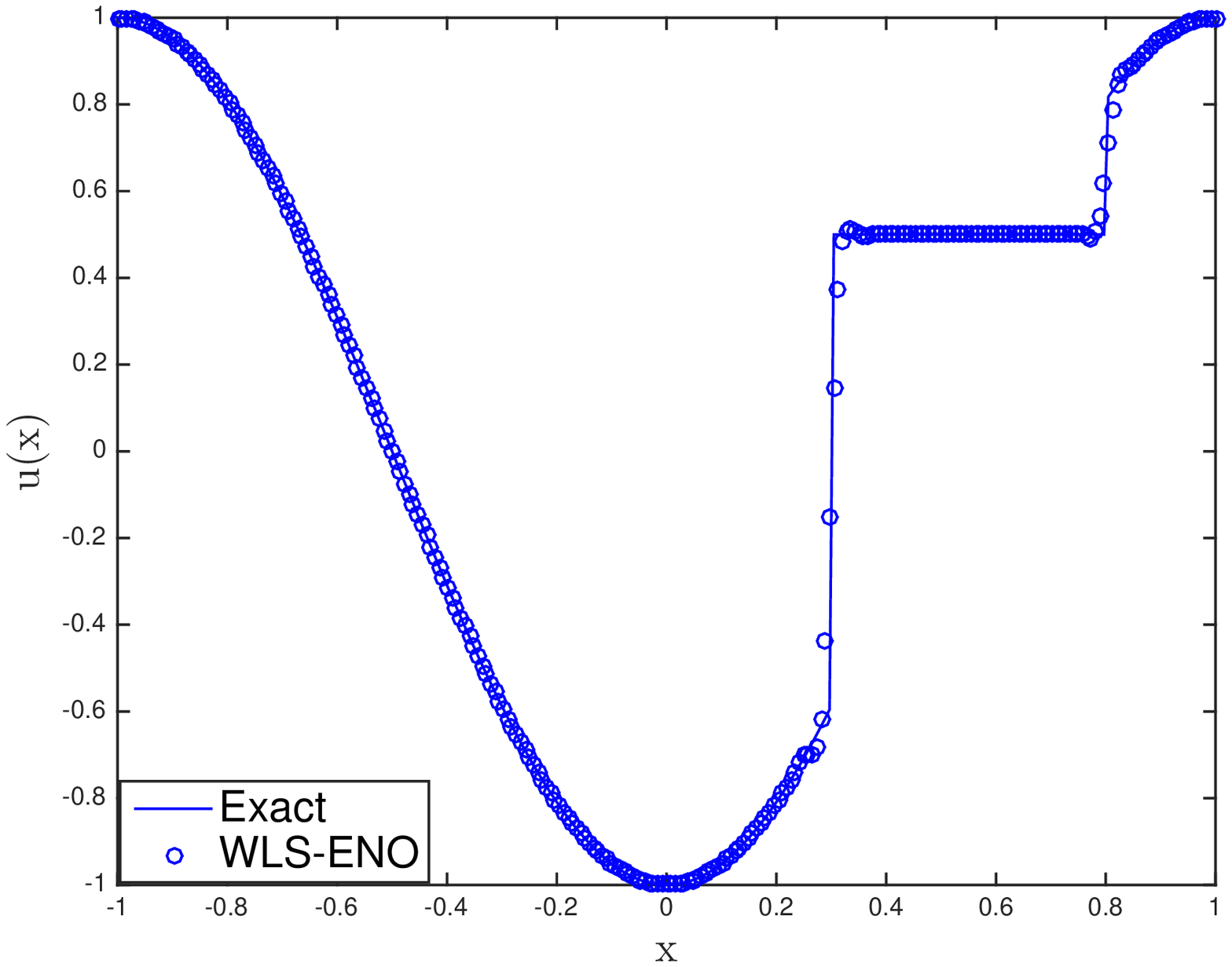}\\
{\small{}(b) Numerical solution at $t=0.5$. }
\par\end{center}{\small \par}%
\end{minipage}
\par\end{centering}

\centering{}\protect\caption{\label{fig:discont-ic-waveeq}Discontinuous initial condition (left)
and numerical solution with WLS-ENO (right) at $t=0.5$ for the linear
wave equation.}
\end{figure}

To demonstrate the accuracy and stability of WLS-ENO for discontinuous
solutions, we change the initial condition to be a piecewise smooth
function
\begin{equation}
u(x,0)=\begin{cases}
\sin(\pi x) & -1\leq x<-0.2\cup0.3<x\leq1\\
0.5 & -0.2\leq x\leq0.3
\end{cases},
\end{equation}
as shown in Figure~\ref{fig:discont-ic-waveeq}(a). Figure~\ref{fig:discont-ic-waveeq}(b)
shows the solution at $t=0.5$ using WLS-ENO. The results show that
WLS-ENO scheme well preserved the sharp feature and maintained the
non-oscillatory property.

\subsubsection{1-D Burgers' Equation}

\begin{figure}
\begin{centering}
\begin{minipage}[t]{0.49\columnwidth}%
\begin{center}
\includegraphics[width=1\textwidth]{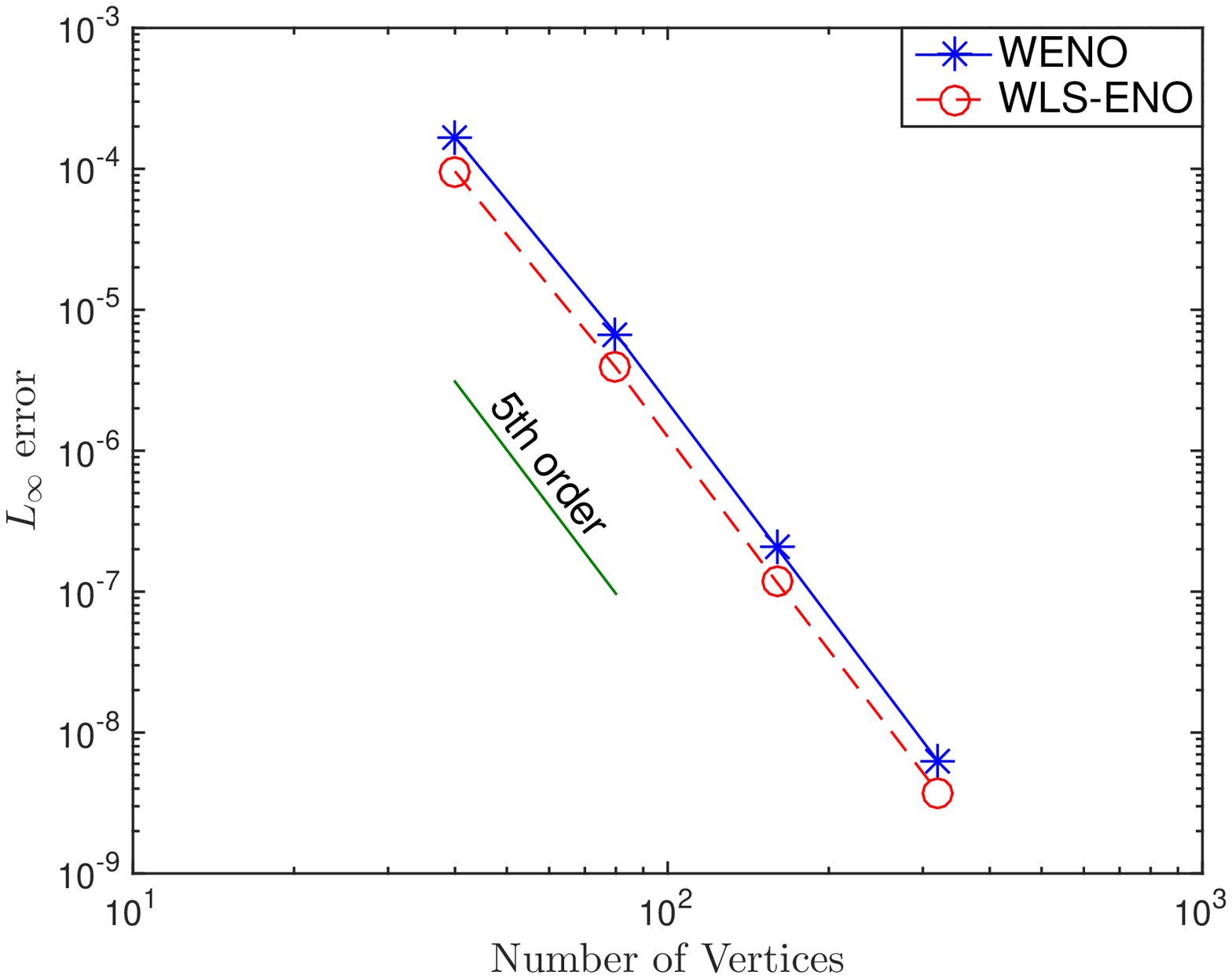}\\
{\small{}(a)}\label{fig:burger1D}{\small{} Errors on uniform grids
at $t=1$.}
\par\end{center}{\small \par}%
\end{minipage}\hfill %
\begin{minipage}[t]{0.49\columnwidth}%
\begin{center}
\includegraphics[width=1\textwidth]{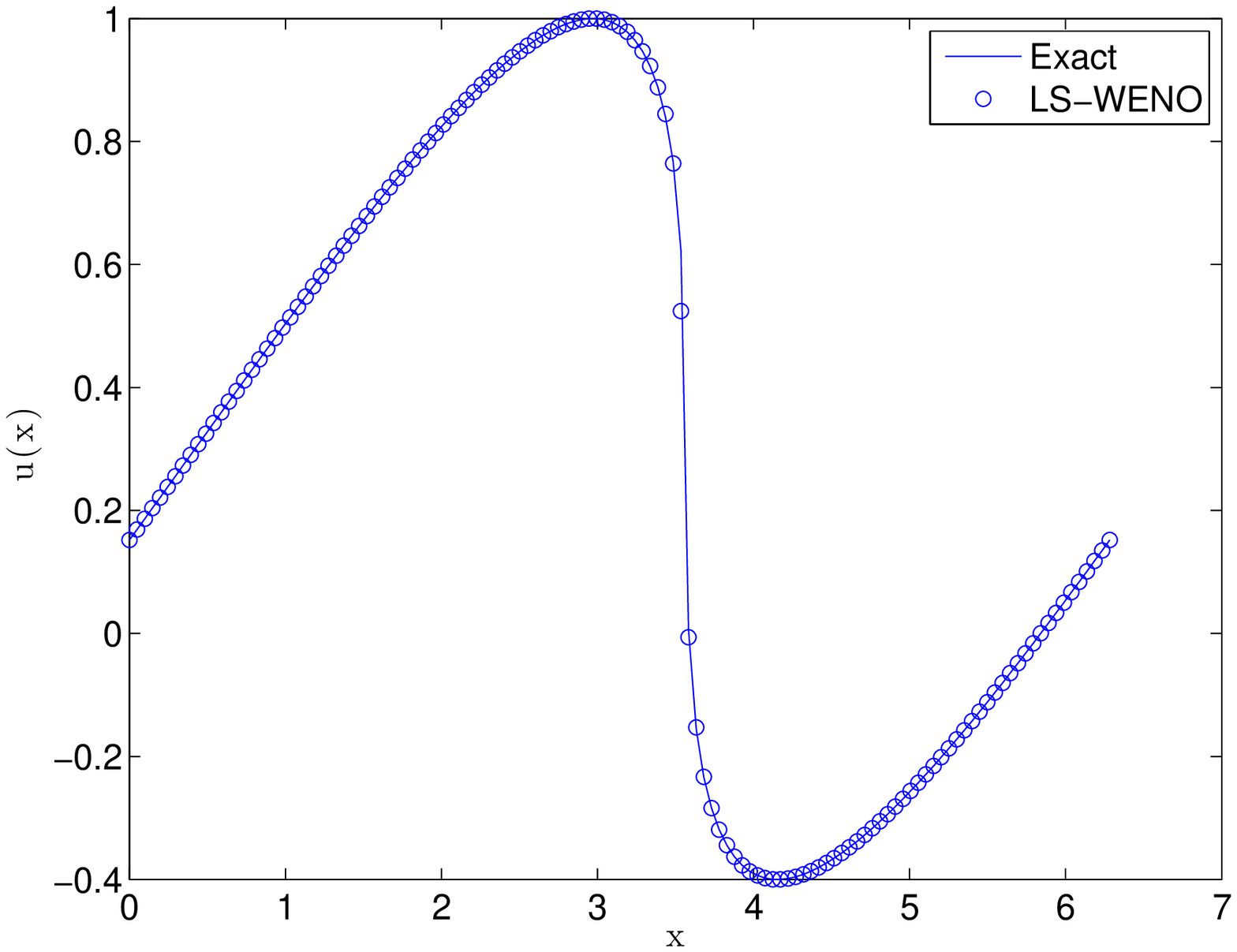}\\
{\small{}(b) }\label{fig:burger1D_breaking}{\small{}Numerical and
exact solution at $t=1.4$. }
\par\end{center}{\small \par}%
\end{minipage}
\par\end{centering}

\centering{}\protect\caption{\label{fig:Comparison-of-WENO1D}Comparison of WENO and WLS-ENO scheme
for 1-D Burgers' equation at $t=1$ (left) and numerical solution
with WLS-ENO (right) at $t=1.4$ }
\end{figure}

Next, we test WLS-ENO with the 1-D Burgers' equation,
\begin{equation}
\frac{\partial u}{\partial t}+u\frac{\partial u}{\partial x}=0,\:0\leq x\leq2\pi,
\end{equation}
with periodic boundary conditions and the initial condition
\begin{equation}
u(x,0)=0.3+0.7\sin(x),\:0\leq x\leq2\pi.
\end{equation}
Although the initial condition is smooth, a discontinuity develops
at time $t=1.4$. To assess the order of accuracy, Figure~\ref{fig:Comparison-of-WENO1D}(a)
shows the solutions from fifth-order WLS-ENO, compared with fifth-order
WENO under grid refinement at time $t=1$, starting from a grid with
64 grid points. The results show that both WLS-ENO and WENO delivered
fifth order accuracy, while WLS-ENO is more accurate. Figure~\ref{fig:Comparison-of-WENO1D}(b)
shows the numerical solution from WLS-ENO overlaid on top of the exact
solution at $t=1.4$. We can see that WLS-ENO scheme approximated
the solution very well.

\subsubsection{1-D Euler Equations: Sod's Problem}

The above tests demonstrate the accuracy and stability of WLS-ENO
for 1-D benchmark problems. For a more realistic problem, we consider
the 1-D Euler equation 
\begin{equation}
\left(\begin{array}{c}
\rho\\
\rho v\\
E
\end{array}\right)_{t}+\left(\begin{array}{c}
\rho v\\
\rho v^{2}+p\\
v\left(E+p\right)
\end{array}\right)_{x}=\vec{0},\label{eq:EulerEq}
\end{equation}
with the equation of state for ideal polytropic gas
\begin{equation}
E=\frac{p}{\gamma-1}+\frac{1}{2}\rho v^{2},
\end{equation}
where $\rho$ denotes the gas density, $v$ the velocity, $p$ the
pressure, $E$ the energy, and $\gamma=1.4$ a constant specific to
air. We perform characteristic decomposition for the Euler equation
\cite{Leveque98nonlinearconservation} and solve the conservation
law characteristic-wise using the fifth order WLS-ENO scheme on an
unstructured (i.e., non-uniform) grid. In more detail, if we introduce
the speed of sound $c$ by 
\begin{equation}
c=\sqrt{\frac{\gamma p}{\rho}},
\end{equation}
and enthalpy $H$ by
\begin{equation}
H=\frac{E+p}{\rho},
\end{equation}
we have the eigenvalue decomposition for the Jacobian as
\begin{equation}
R^{-1}(u)f'(u)R(u)=\varLambda(u),
\end{equation}
where
\begin{equation}
f'(u)=\left(\begin{array}{ccc}
0 & 1 & 0\\
\left(\frac{\gamma-3}{2}\right)v^{2} & (3-\gamma)v & \gamma-1\\
\left(\frac{\gamma-1}{2}\right)v^{3}-vH & H-(\gamma-1)v^{2} & \gamma v
\end{array}\right),
\end{equation}
\begin{equation}
\varLambda(u)=\left(\begin{array}{ccc}
v-c\\
 & v\\
 &  & v+c
\end{array}\right),
\end{equation}
and
\begin{equation}
R(u)=\left(\begin{array}{ccc}
1 & 1 & 1\\
v-c & v & v+c\\
H-cv & \frac{1}{2}v^{2} & H+cv
\end{array}\right).
\end{equation}

\begin{figure}
\centering{}%
\begin{minipage}[t]{0.33\columnwidth}%
\begin{center}
\includegraphics[width=1\textwidth]{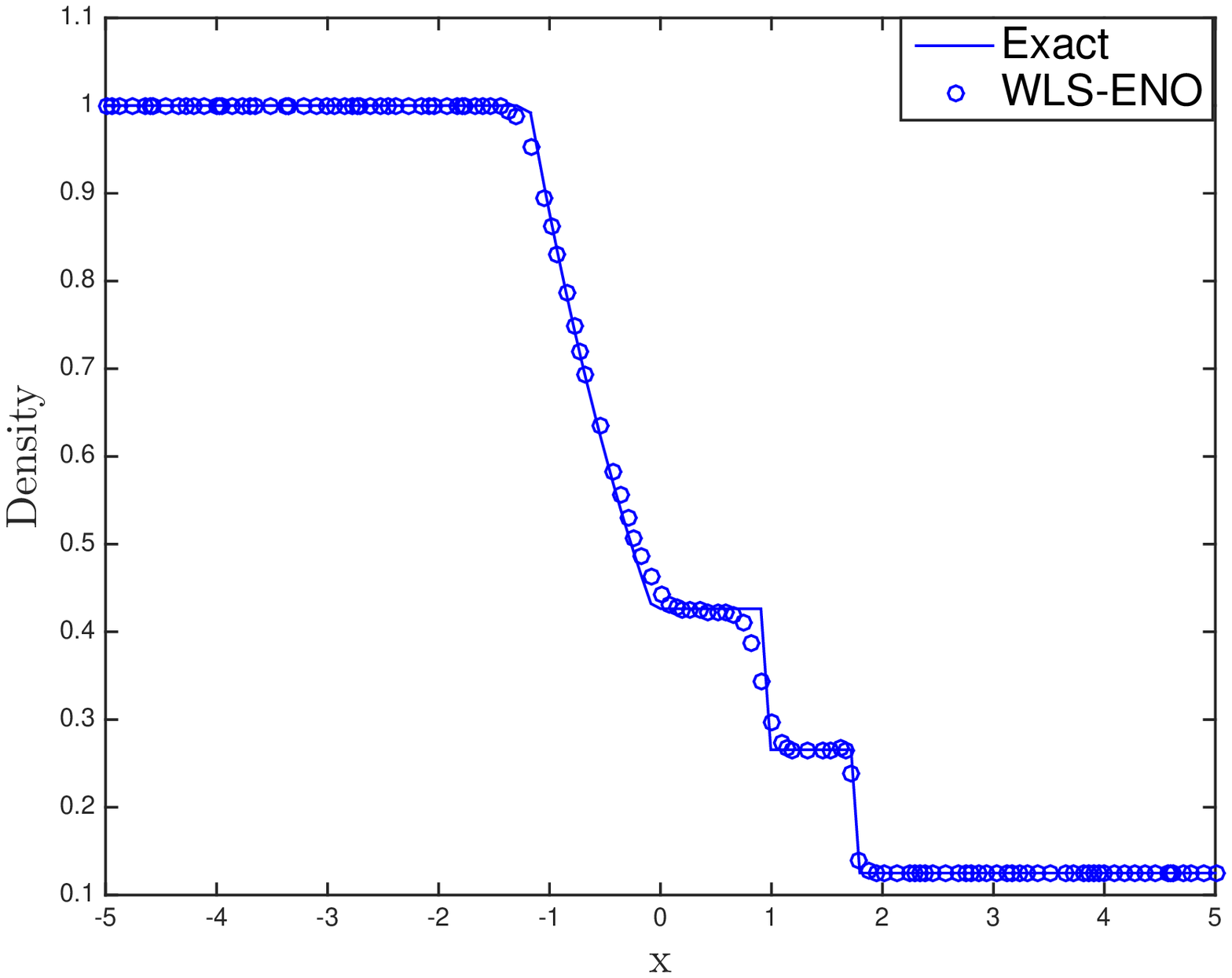}\\
{\small{}(a) Density.}
\par\end{center}{\small \par}%
\end{minipage}%
\begin{minipage}[t]{0.33\columnwidth}%
\begin{center}
\includegraphics[width=1\textwidth]{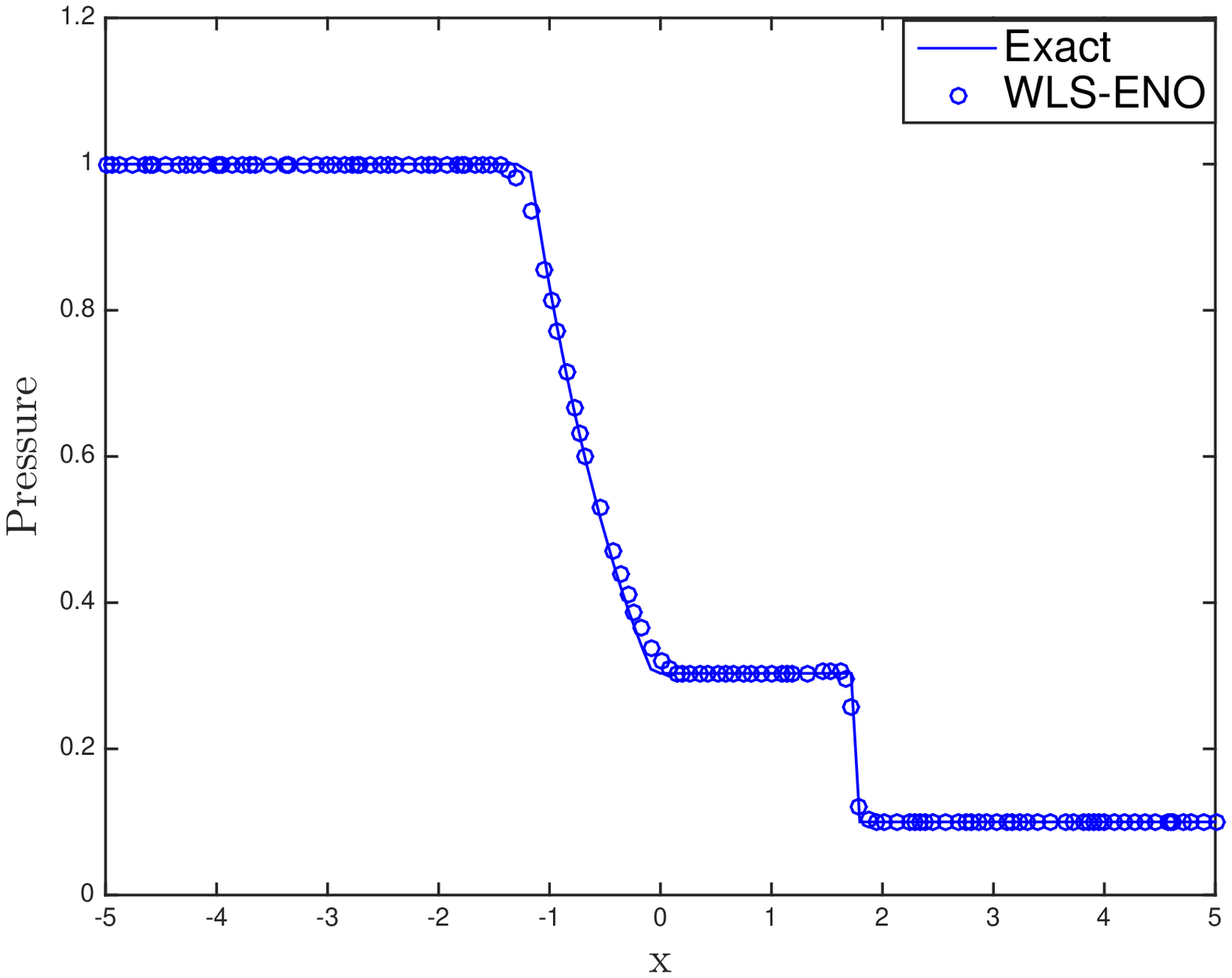}\\
{\small{}(b) Pressure.}
\par\end{center}{\small \par}%
\end{minipage}%
\begin{minipage}[t]{0.33\columnwidth}%
\begin{center}
\includegraphics[width=1\textwidth]{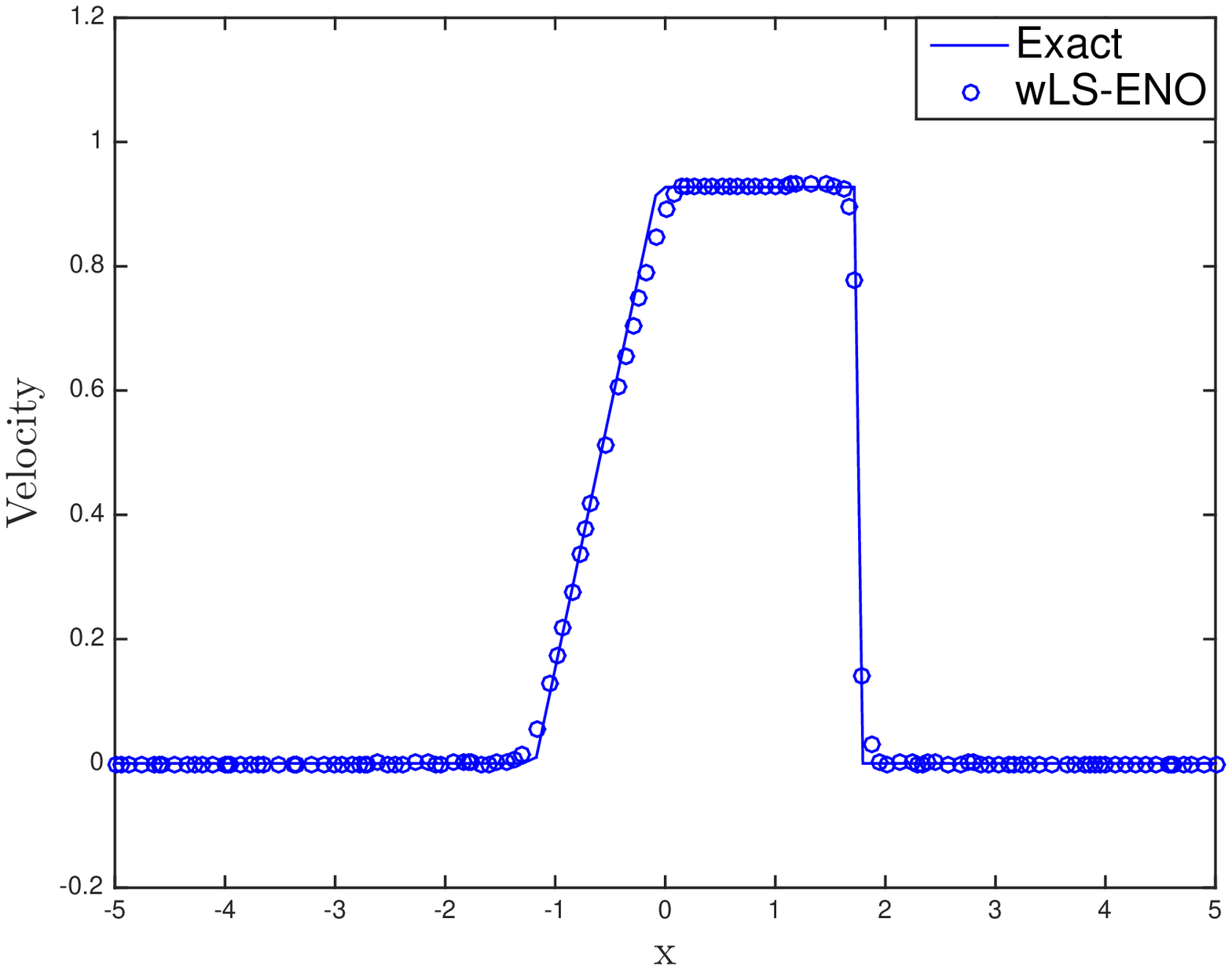}\\
{\small{}(c) Velocity.}
\par\end{center}{\small \par}%
\end{minipage}\protect\caption{\label{fig:Solutions-of-Euler}Solutions of 1-D Euler equations (\ref{eq:EulerEq})
at $t=1$ by fifth order WLS-ENO on non-uniform grid.}
\end{figure}

We compute the density, velocity and pressure with the initial condition
given by Sod's problem \cite{ShuWeno98}
\begin{equation}
\left(\rho_{L},v_{L},p_{L}\right)=\left(1,0,1\right),\left(\rho_{R},v_{R},p_{R}\right)=\left(0.125,0,0.1\right)
\end{equation}
Figure~\ref{fig:Solutions-of-Euler} shows the numerical solutions
at state $t=1$ compared against the exact solution. From the results,
we can see that the results match the exact solution very well and
are non-oscillatory near discontinuities.

\subsubsection{Interacting Blast Waves}

\begin{figure}
\centering{}%
\begin{minipage}[t]{0.33\columnwidth}%
\begin{center}
\includegraphics[width=1\textwidth]{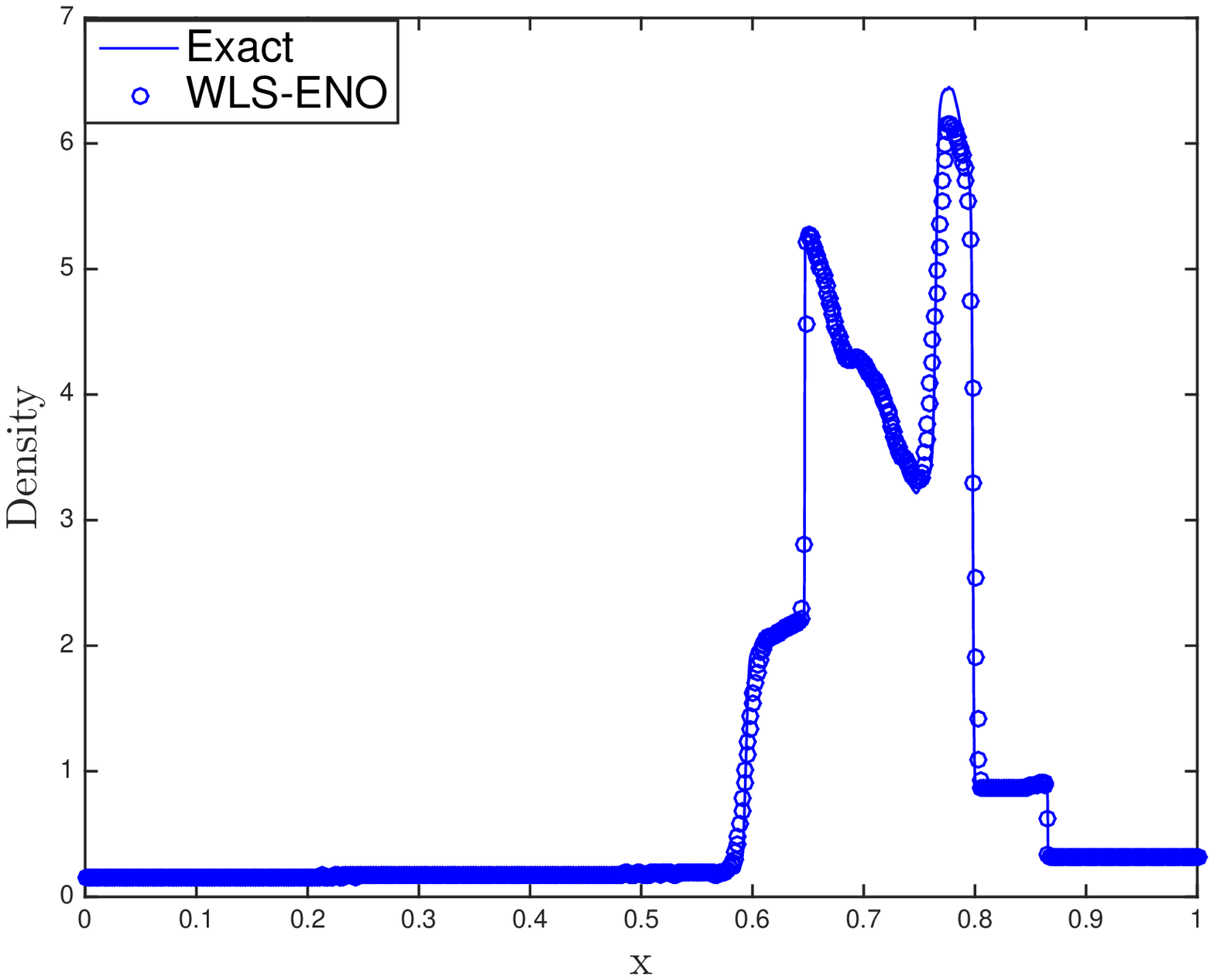}\\
{\small{}(a) Density.}
\par\end{center}{\small \par}%
\end{minipage}%
\begin{minipage}[t]{0.33\columnwidth}%
\begin{center}
\includegraphics[width=1\textwidth]{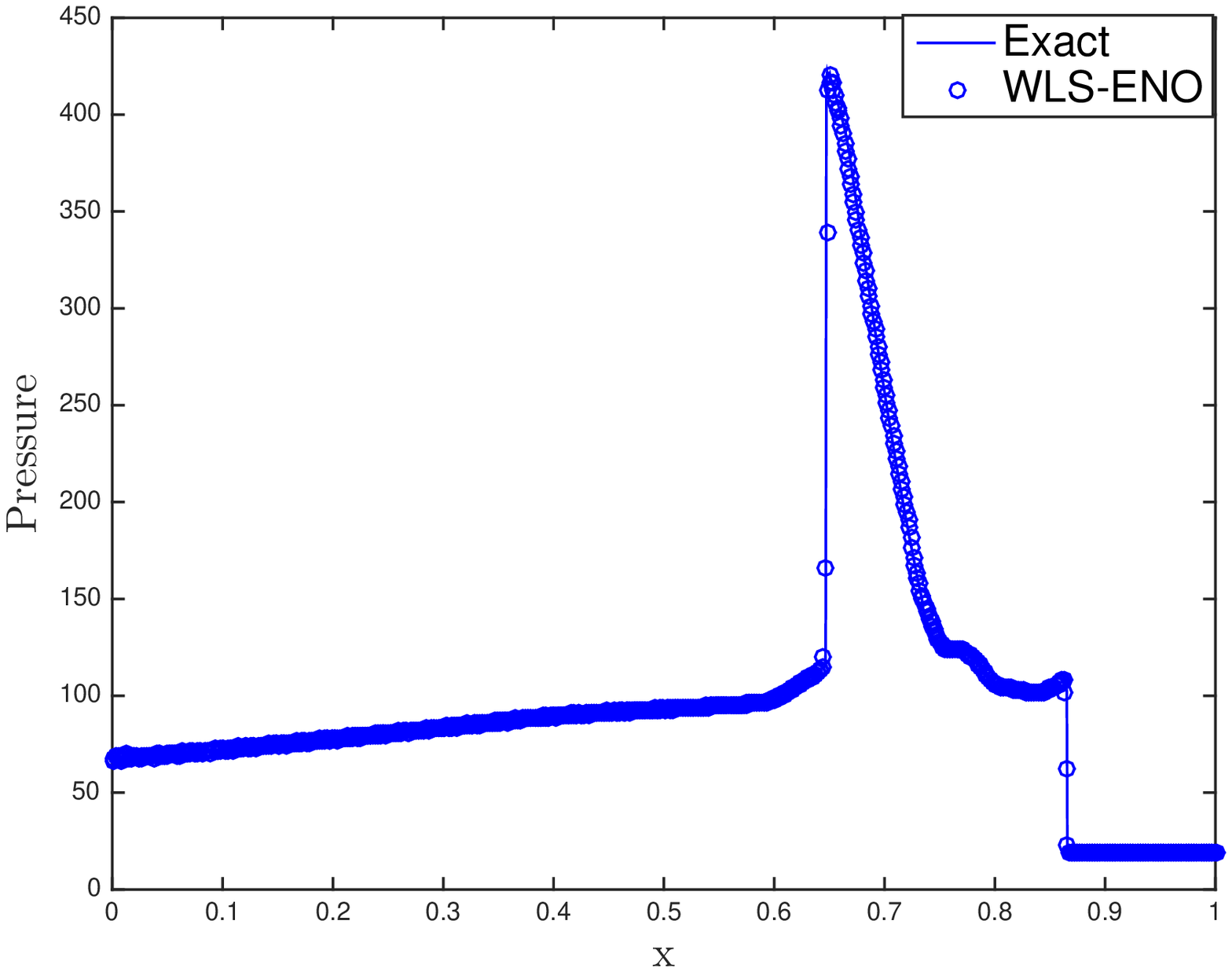}\\
{\small{}(b) Pressure.}
\par\end{center}{\small \par}%
\end{minipage}%
\begin{minipage}[t]{0.33\columnwidth}%
\begin{center}
\includegraphics[width=1\textwidth]{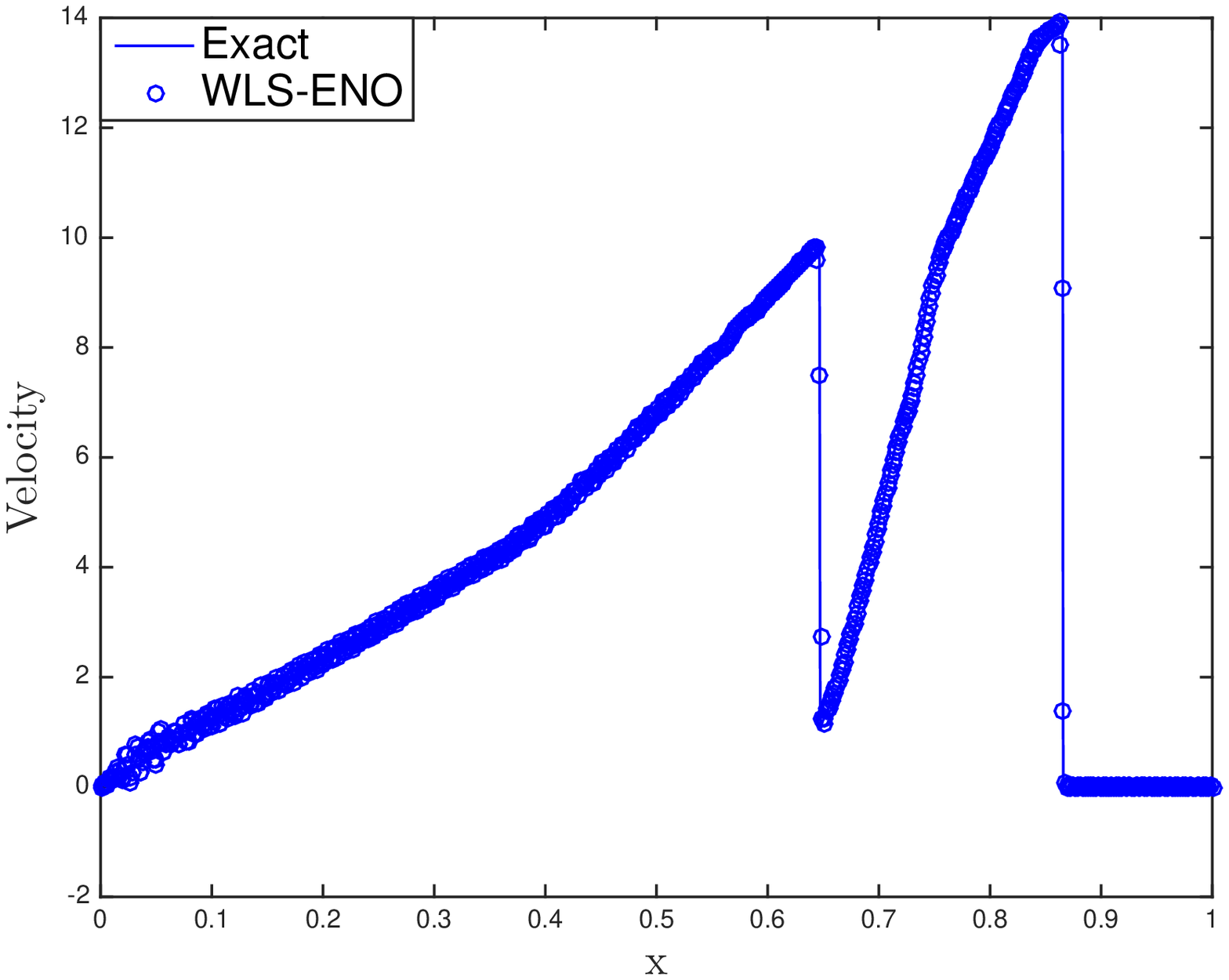}\\
{\small{}(c) Velocity.}
\par\end{center}{\small \par}%
\end{minipage}\protect\caption{\label{fig:interactive blasting wave}Solutions of 1-D interacting
blast waves at $t=0.038$ by fifth order WLS-ENO. }
\end{figure}

This is a test problem specially designed by Woodward and Collela
to illustrate the strong relationship between the accuracy of the
overall flow solution and the sharp resolution of discontinuities.
It has the following the initial condition:
\begin{equation}
\left(\rho,u,P\right)=\begin{cases}
\begin{array}{cc}
\left(1,0,1000\right) & 0\leq x<0.1\\
\left(1,0,0.01\right) & 0.1\leq x<0.9.\\
\left(1,0,100\right) & 0.9\leq x<1
\end{array}\end{cases}
\end{equation}
Reflective boundary conditions are applied at both sides. For this
example, we use fifth order WLS-ENO scheme and plot the results at
time $t=0.038$. Figure~\ref{fig:interactive blasting wave} shows
excellent agreement of the numerical solution with the exact solution.

\subsection{2-D Results}

We now present results of WLS-ENO with 2-D unstructured meshes for
problems with smooth or piecewise smooth solutions, including the
wave equation, Burgers' equation, and the Euler equations with two
different initial conditions.

\subsubsection{2-D Wave Equation}

\begin{figure}
\begin{minipage}[t]{0.49\columnwidth}%
\begin{center}
\includegraphics[width=1\textwidth]{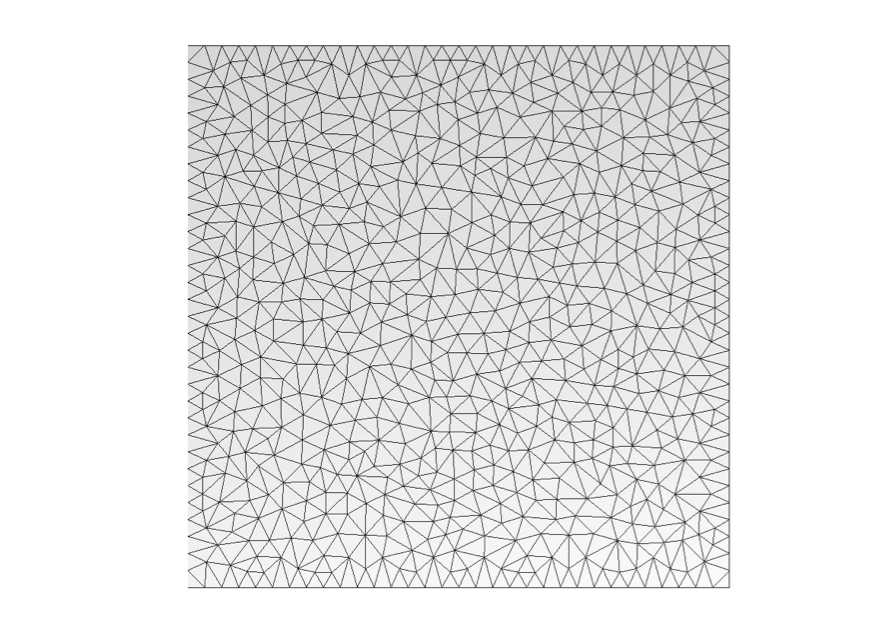}
\par\end{center}

\begin{center}
\protect\caption{\label{fig:sample_2d_mesh}A sample unstructured triangular mesh for
solving 2-D wave equation.}

\par\end{center}%
\end{minipage}\hfill%
\begin{minipage}[t]{0.49\columnwidth}%
\begin{center}
\includegraphics[width=1\textwidth]{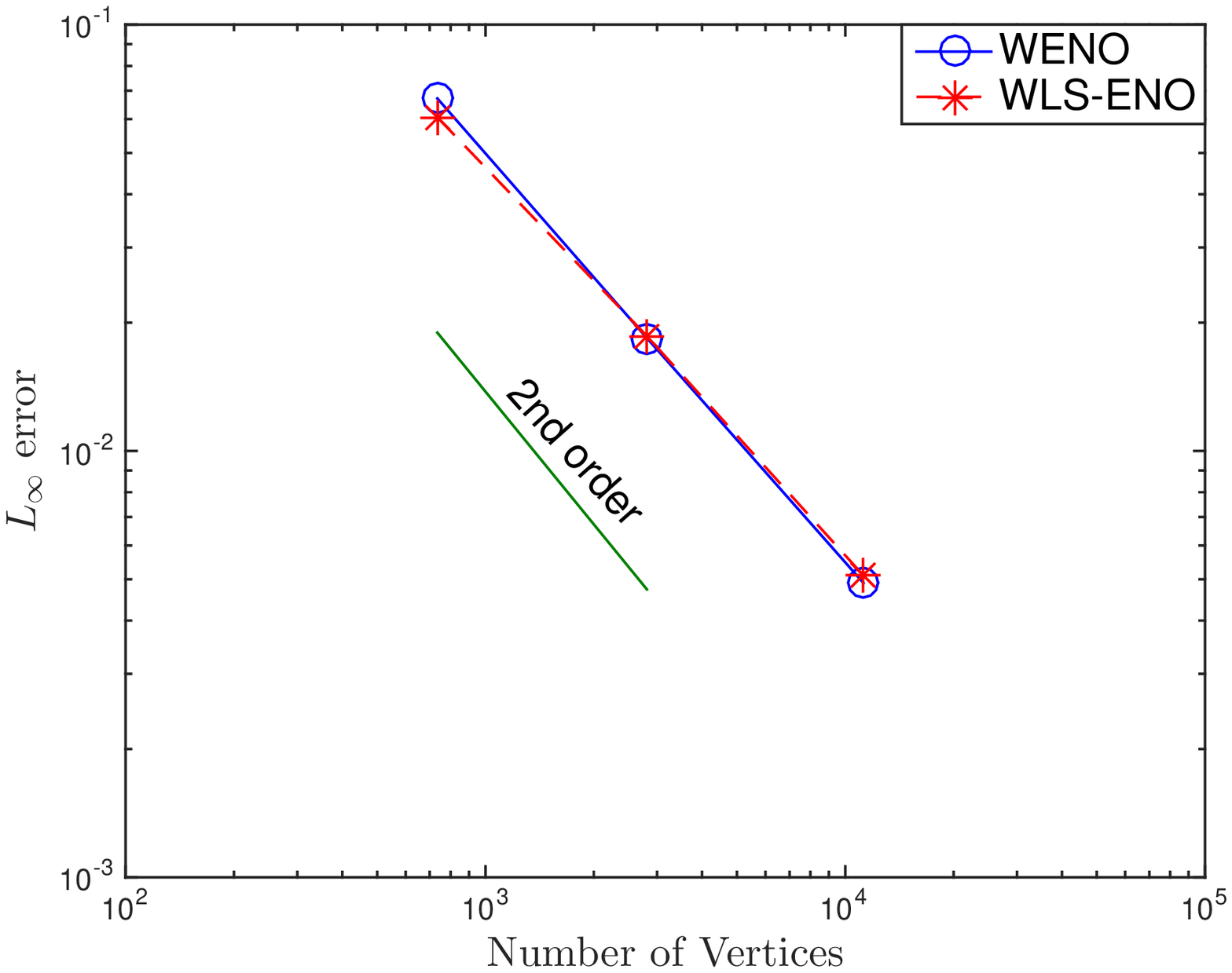}
\par\end{center}

\begin{center}
\protect\caption{\label{fig:Errors-in-solutions2D}Errors of solutions of 2-D wave
equation with WLS-ENO at $t=1$.}

\par\end{center}%
\end{minipage}
\end{figure}

\begin{figure}
\centering{}%
\begin{minipage}[t]{0.5\columnwidth}%
\begin{center}
\includegraphics[width=1\textwidth]{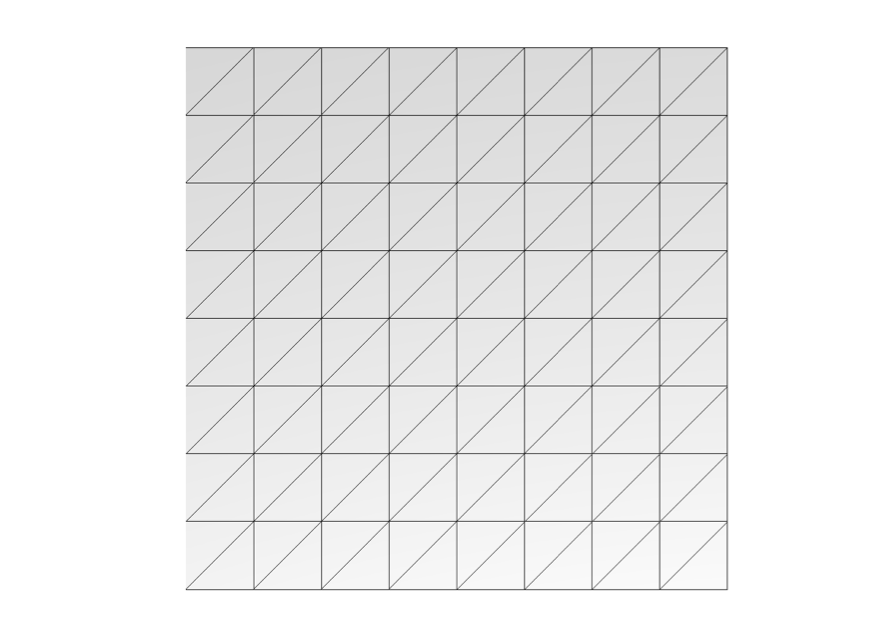}
\par\end{center}

\begin{center}
\protect\caption{\label{fig:uniform-trimesh}Sample uniform triangular mesh for solving
2-D wave equation.}

\par\end{center}%
\end{minipage}\hfill%
\begin{minipage}[t]{0.49\columnwidth}%
\begin{center}
\includegraphics[width=1\textwidth]{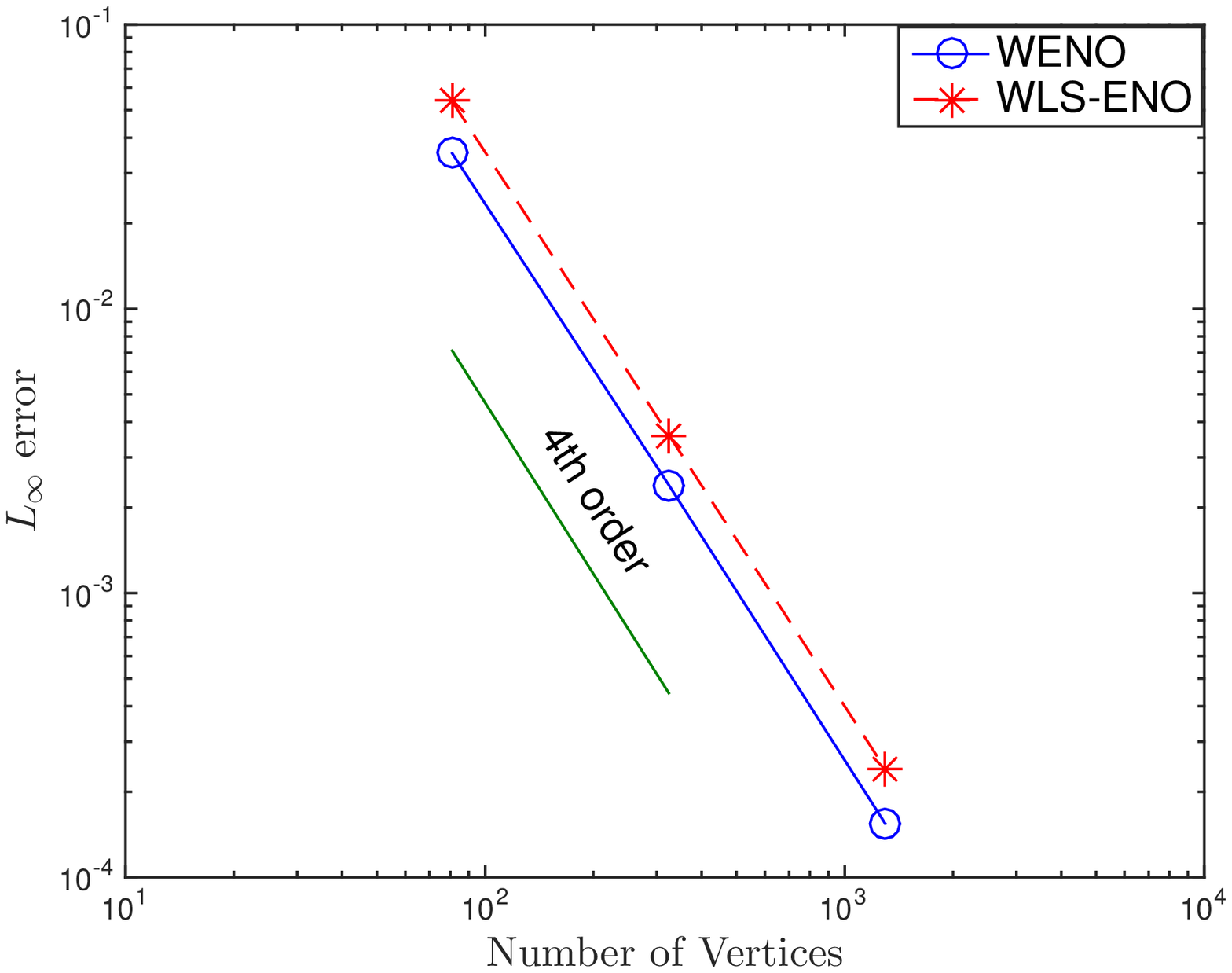}
\par\end{center}

\begin{center}
\protect\caption{\label{fig:Errors-of-wave}Errors of 2-D wave equation of WENO and
WLS-ENO on uniform meshes.}

\par\end{center}%
\end{minipage}
\end{figure}

As in 1-D, we first consider the wave equation,
\begin{equation}
u_{t}+u_{x}+u_{y}=0,
\end{equation}
with periodic boundary conditions and the initial condition
\begin{equation}
u_{0}(x,y)=\sin\left(\frac{\pi}{2}\left(x+y\right)\right),\qquad-2\leq x\leq2,\:-2\leq y\leq2.
\end{equation}
For this problem, the solution remains smooth over time. We solve
the problem using third-order WLS-ENO scheme and third-order WENO
scheme. Figure~\ref{fig:Errors-in-solutions2D} shows the errors
at $t=1$, and it can be seen that both methods achieved only second
order convergence. This convergence rate is expected, because the
derivatives can only be approximated to second-order accurate by polynomial
approximations over nonuniform unstructured meshes without symmetry.
When applying WLS-ENO on a uniform mesh, such as that shown in Figure~\ref{fig:uniform-trimesh},
it would deliver the convergence rate one order higher due to error
cancellation, similar to WENO and other finite difference methods,
as illustrated with the fourth-order WLS-ENO scheme and fourth-order
WENO scheme in Figure~\ref{fig:Errors-of-wave}. Note that on uniform
meshes, WLS-ENO may be slightly less accurate than WENO because it
uses a larger stencil.

\subsubsection{2-D Burgers' Equation}

\begin{figure}
\begin{centering}
\begin{minipage}[t]{0.49\columnwidth}%
\begin{center}
\includegraphics[width=1\textwidth]{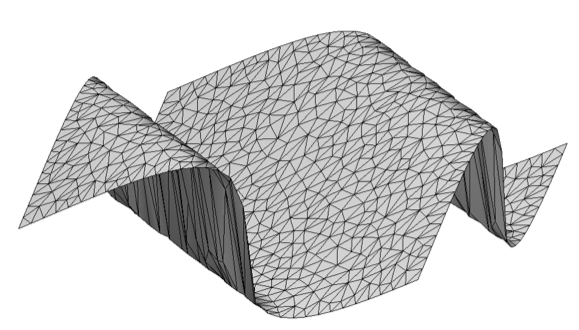}
\par\end{center}%
\end{minipage}\hfill%
\begin{minipage}[t]{0.49\columnwidth}%
\begin{center}
\includegraphics[width=1\textwidth]{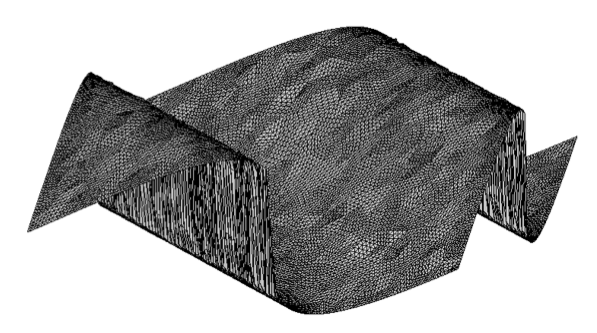}
\par\end{center}%
\end{minipage}
\par\end{centering}

\centering{}\protect\caption{\label{fig:Burgers2-D}Exact solution (left) and numerical solution
(right) with fourth-order WLS-ENO for the 2-D Burgers' equation at
$t=0.5$. }
\end{figure}

For piecewise smooth solutions, we solve the 2-D Burgers' equation
\begin{equation}
u_{t}+\left(\frac{u^{2}}{2}\right)_{x}+\left(\frac{u^{2}}{2}\right)_{y}=0
\end{equation}
over $[-2,2]^{2}$, with periodic boundary conditions and the initial
condition 
\begin{equation}
u_{0}(x,y)=0.3+0.7\sin\left(\frac{\pi}{2}\left(x+y\right)\right).
\end{equation}
Although the initial condition is smooth, discontinuities develop
over time. Figure~\ref{fig:Burgers2-D}(left) shows the exact solution
at $t=0.5$, when the solution becomes discontinuous, and Figure~\ref{fig:Burgers2-D}(right)
shows the result of fourth-order WLS-ENO scheme under non-uniform
grid refinement. This result shows that the overall solution maintains
accurate as discontinuities develop. 

\subsubsection{2-D Vortex Evolution Problem}

This is one of the few problems that has exact solutions for the compressible
Euler equations. The test case involves the convection of an isentropic
vortex in inviscid flow. It is used to show the ability of numerical
schemes to accurately capture vortical flows. The 2-D Euler equations
have the following form
\begin{equation}
\left(\begin{array}{c}
\rho\\
\rho u\\
\rho v\\
E
\end{array}\right)_{t}+\left(\begin{array}{c}
\rho u\\
\rho u^{2}+p\\
\rho uv\\
u\left(E+p\right)
\end{array}\right)_{x}+\left(\begin{array}{c}
\rho v\\
\rho uv\\
\rho v^{2}+p\\
v\left(E+p\right)
\end{array}\right)_{y}=\vec{0},
\end{equation}
where
\begin{equation}
E=\frac{p}{\gamma-1}+\frac{1}{2}\rho\left(u^{2}+v^{2}\right).
\end{equation}

\begin{figure}
\begin{centering}
\includegraphics[width=0.5\textwidth]{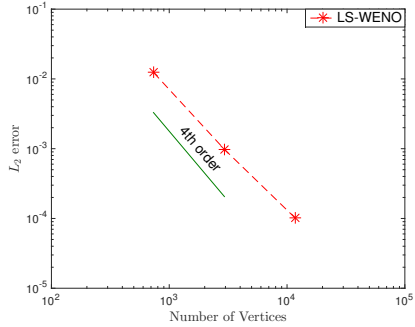}
\par\end{centering}

\centering{}\protect\caption{\label{fig:Density-2D Vortex Evolution}Density errors in 2-D Vortex
Evolution problem by fifth order WLS-ENO on triangular meshes.}
\end{figure}

We consider an idealized problem: The mean flow is $\rho_{\infty}=1$,
$p_{\infty}=1$, $(u_{\infty},v_{\infty})=(1,1)$. As an initial condition,
an isentropic vortex with no perturbation in entropy ($\delta S=0$)
is added to the mean flow field. The perturbation values are given
by
\begin{equation}
\left(\delta u,\delta v\right)=\frac{\beta}{2\pi}e^{\frac{1-r^{2}}{2}}\left(-\overline{y},\overline{x}\right),
\end{equation}
\begin{equation}
\delta T=-\frac{\left(\gamma-1\right)\beta^{2}}{8\gamma\pi^{2}}e^{1-r^{2}}
\end{equation}
where $\left(\overline{x},\overline{y}\right)=\left(x-5,y-5\right)$,
$r^{2}=\overline{x}^{2}+\overline{y}^{2}$, and the vortex strength
$\beta=5$. The computational domain is taken as $\left[0,10\right]\times\left[0,10\right]$,
extended periodically in both directions. The exact solution of this
problem is simply a vortex convecting with the speed $(1,1)$ in the
diagonal direction. In this example, we use fifth-order WLS-ENO scheme
and compute the results up to $t=1$. Figure~\ref{fig:Density-2D Vortex Evolution}
 shows that the result achieves the expected convergence rate.

\subsubsection{2-D Explosion Test Problem}

In this test, we solve the Euler equation with the initial condition
\begin{equation}
\left(\rho,u,v,p\right)^{T}=\begin{cases}
\left(1,0,0,1\right)^{T} & \sqrt{x^{2}+y^{2}}\leq0.2\\
\left(0.125,0,0.1\right)^{T} & \sqrt{x^{2}+y^{2}}>0.2
\end{cases}.
\end{equation}
We solve the problem on a domain of a unit disk centered at the origin,
triangulated with meshes similar to that in Figure~\ref{fig:sample_2d_mesh}.
For this test, we utilize characteristic decomposition for the Jacobian
of Euler system \cite{HuShu99WENO} and run the test up to $t=0.1$
to ensure that the explosion wave do not reach the boundary. 

\begin{figure}
\centering{}%
\begin{minipage}[t]{0.49\columnwidth}%
\begin{center}
\includegraphics[width=1\textwidth]{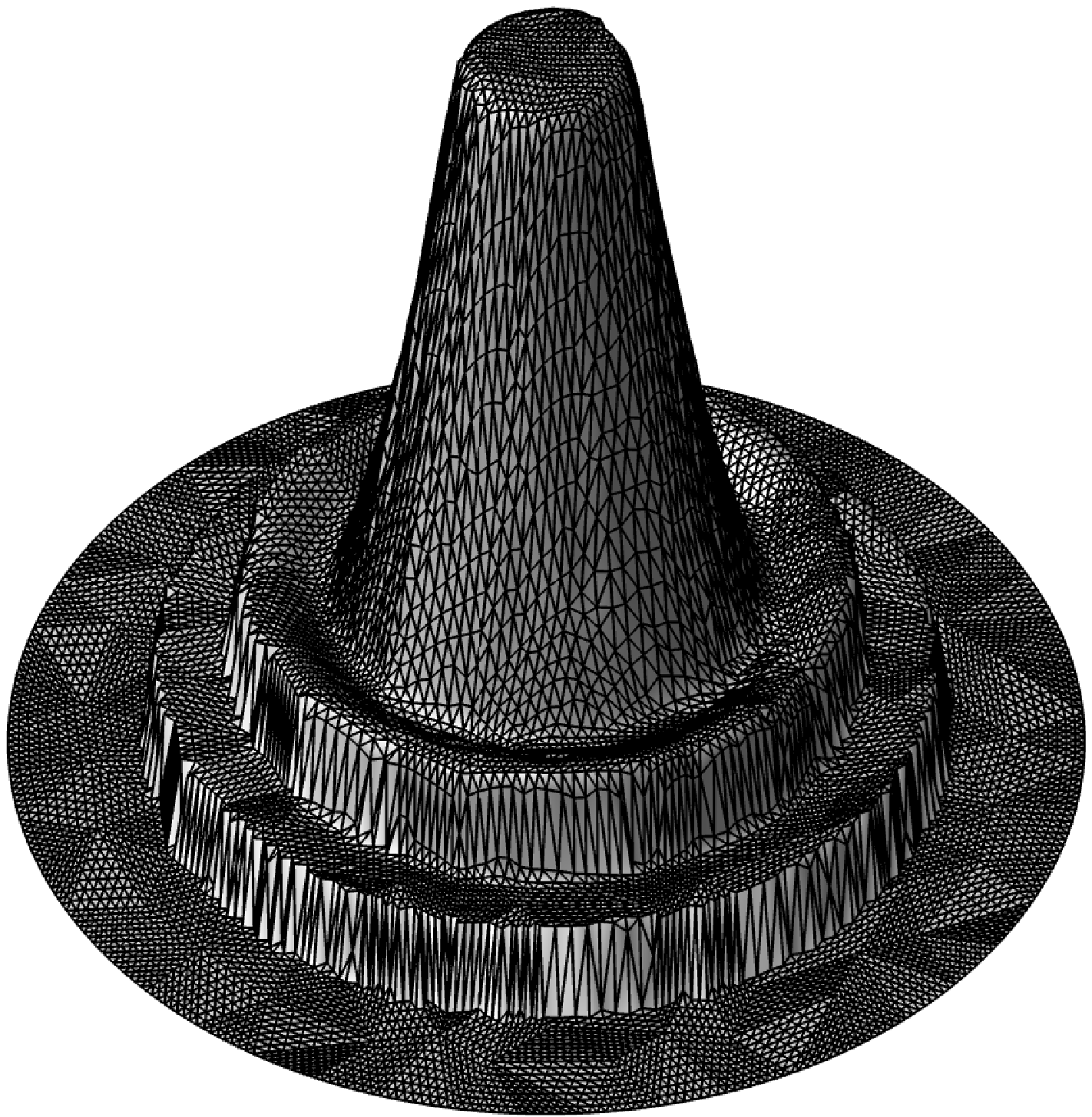}
\par\end{center}

\begin{center}
\protect\caption{\label{fig:2D Explosion}Numerical solution of 2-D explosion test
by third-order WLS-ENO at $t=0.1$. }

\par\end{center}%
\end{minipage}\hfill%
\begin{minipage}[t]{0.49\columnwidth}%
\begin{center}
\includegraphics[width=1\textwidth]{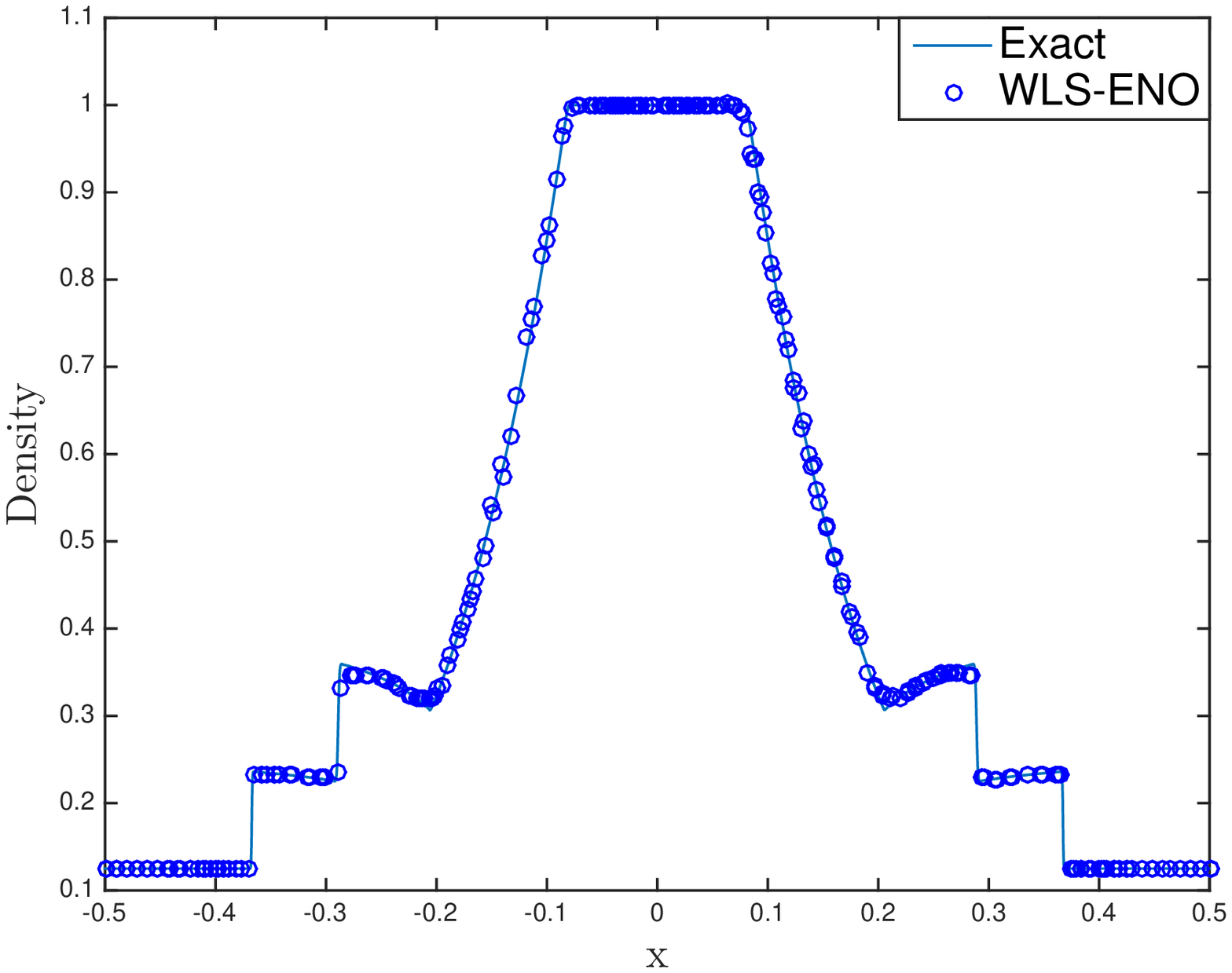}
\par\end{center}

\begin{center}
\protect\caption{\label{fig:2D Explosion-1D compare}Numerical solution along $x$
axis of 2-D explosion test by third-order WLS-ENO versus analytic
solution at $t=0.1$.}

\par\end{center}%
\end{minipage}
\end{figure}

Because of symmetry, this example is equivalent to the spherical one
dimensional Euler equation with source terms \cite{EFToro},
\begin{equation}
\frac{\partial}{\partial t}\left(\begin{array}{c}
\rho\\
\rho u\\
E
\end{array}\right)+\frac{\partial}{\partial r}\left(\begin{array}{c}
\rho u\\
\rho u^{2}+p\\
u\left(E+p\right)
\end{array}\right)=-\frac{d-1}{r}\left(\begin{array}{c}
\rho u\\
\rho u^{2}\\
u\left(E+p\right)
\end{array}\right),\label{eq:Explosion 1D}
\end{equation}
where $r$ is the radial coordinate. Thus, as a reference solution,
we solve this 1-D problem on a very fine mesh composed of 4000 grid
points and compare it with our numerical solution. Figures~\ref{fig:2D Explosion}
and \ref{fig:2D Explosion-1D compare} show that the numerical solutions
agree with the exact solutions very well.

\subsection{3-D Results}

One advantage of WLS-ENO is that it generalizes to 3-D in a straightforward
fashion. We present some numerical results over unstructured meshes
in 3-D, including the wave equation, Burgers' equation, and the Euler
equations.

\subsubsection{3-D Wave Equation}

We first solve the 3-D linear wave equation 
\begin{equation}
u_{t}+u_{x}+u_{y}+u_{z}=0
\end{equation}
over $[-2,2]^{3}$, with periodic boundary conditions and the initial
condition 
\begin{equation}
u(x,y,z,0)=\sin\left(\frac{\pi}{2}\left(x+y+z\right)\right).
\end{equation}

\begin{figure}
\centering{}%
\begin{minipage}[t]{0.49\columnwidth}%
\begin{center}
\includegraphics[width=1\textwidth]{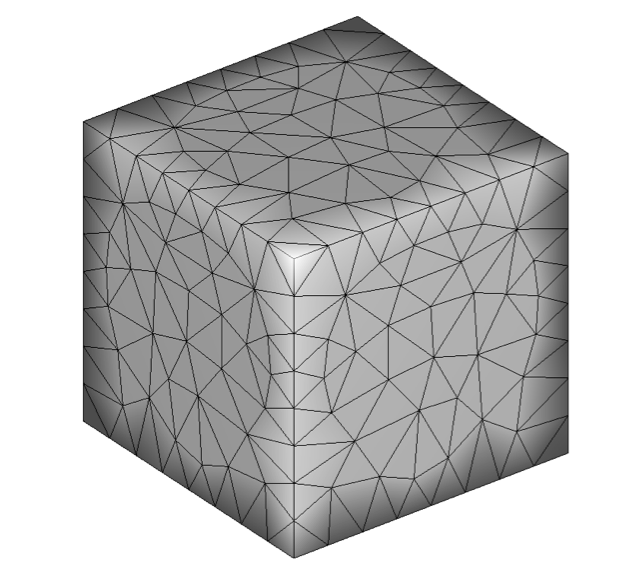}
\par\end{center}

\protect\caption{\label{fig:test-tetmesh}A sample unstructured mesh for solving 3-D
wave equations.}
\end{minipage}\hfill%
\begin{minipage}[t]{0.49\columnwidth}%
\begin{center}
\includegraphics[width=1\textwidth]{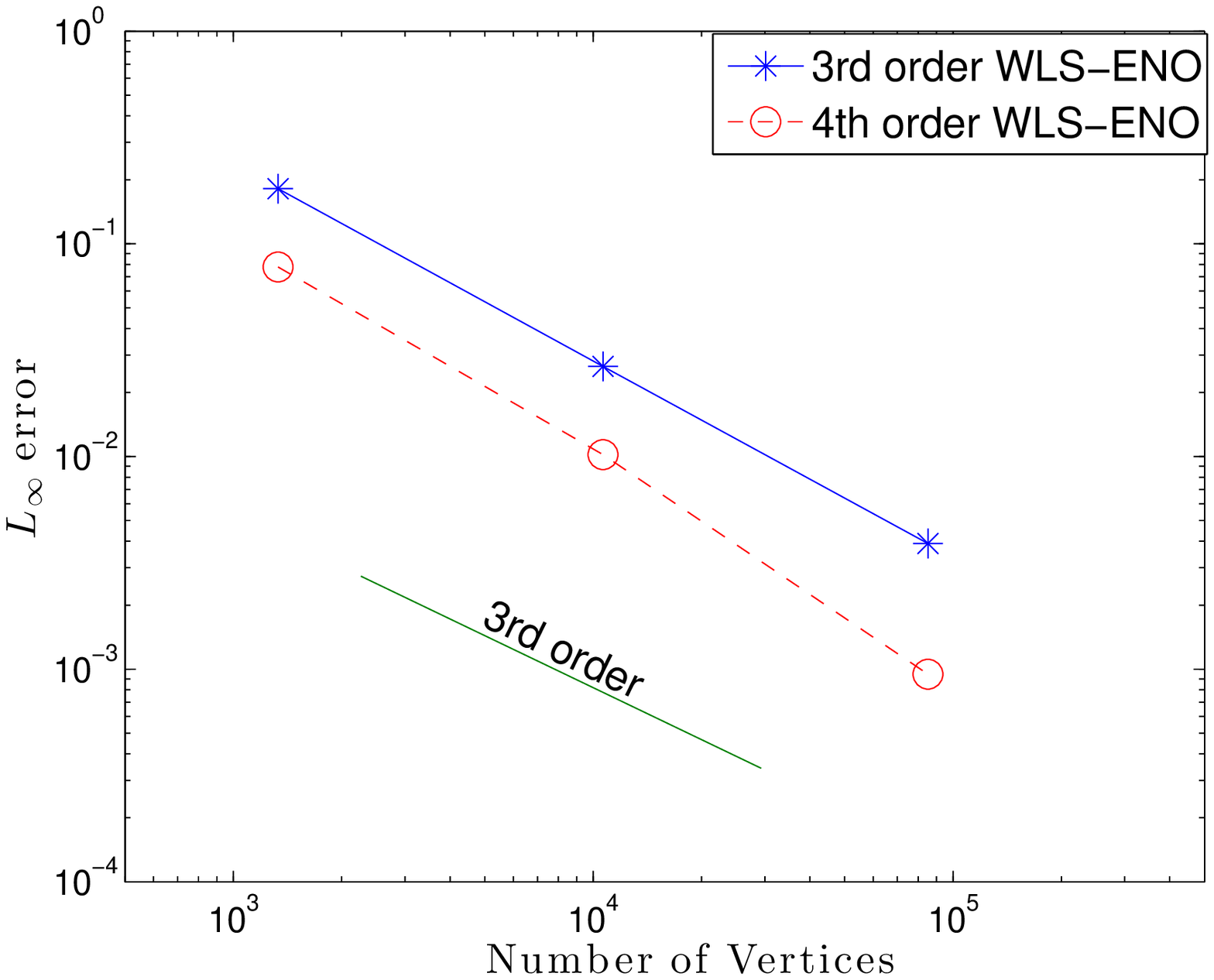}
\par\end{center}

\begin{center}
\protect\caption{\label{fig:conv-3d-wave}Convergence of third- and fourth-order WLS-ENO
schemes for wave equation on tetrahedral meshes.}

\par\end{center}%
\end{minipage}
\end{figure}

We solve the problem using WLS-ENO over a series of unstructured meshes,
where the coarsest mesh is depicted in Figure~\ref{fig:test-tetmesh}.
Figure~\ref{fig:conv-3d-wave} shows the errors with third-order
and fourth-order WLS-ENO schemes under mesh refinement. It is clear
that both schemes achieved the convergence rate close to three, and
the error of fourth order WLS-ENO scheme is about half of that of
the third order scheme. 

\subsubsection{3-D Burgers' Equation}

\begin{figure}[h]
\centering{}%
\begin{minipage}[t]{0.49\columnwidth}%
\begin{center}
\includegraphics[width=1\textwidth]{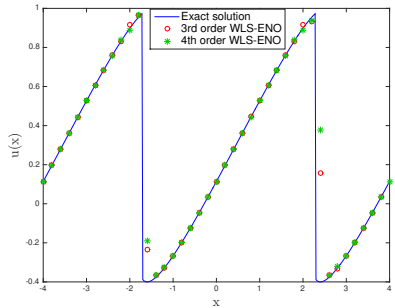}
\par\end{center}

\protect\caption{\label{fig:cross-sections-3d-Burger}1-D cross-sections along $x=y$
and $z=0$ of 3-D Burgers' equation using third- and fourth-order
WLS-ENO schemes at $t=0.5$.}
\end{minipage}\hfill%
\begin{minipage}[t]{0.49\columnwidth}%
\includegraphics[width=1\textwidth]{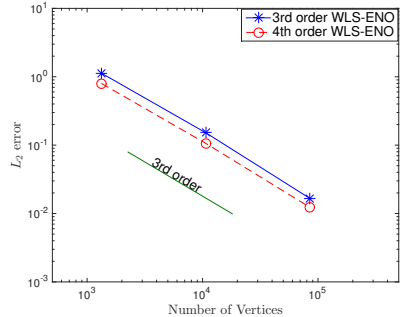}

\begin{center}
\protect\caption{Errors of third- and fourth-order WLS-ENO schemes at $t=0.5$ away
from singularities.}

\par\end{center}%
\end{minipage}
\end{figure}

As a final example, we solve the 3-D nonlinear Burgers' equation
\begin{equation}
u_{t}+\left(\frac{u^{2}}{2}\right)_{x}+\left(\frac{u^{2}}{2}\right)_{y}+\left(\frac{u^{2}}{2}\right)_{z}=0
\end{equation}
over $[-2,2]^{3}$, also with periodic boundary conditions and the
initial condition
\begin{equation}
u\left(x,y,z,0\right)=0.3+0.7\sin\left(\frac{\pi}{2}\left(x+y+z\right)\right).
\end{equation}
Similar to the 2-D case, discontinuities develop at $t=0.5$. Figure~\ref{fig:cross-sections-3d-Burger}
shows a 1-D cross section of the numerical solutions with third-order
and fourth-order WLS-ENO schemes at $t=0.5$, overlaid with the exact
solution. Both solutions are non-oscillatory. In contrast, the third-order
WENO scheme in \cite{zhang2009third} is unstable over a non-uniform
unstructured mesh, and it delivers comparable accuracy to WLS-ENO
over a tetrahedral mesh obtained by decomposing a structured mesh. 

\subsubsection{3-D Explosion Test Problem}

\begin{figure}[H]
\centering{}%
\begin{minipage}[t]{0.49\columnwidth}%
\begin{center}
\includegraphics[width=1\textwidth]{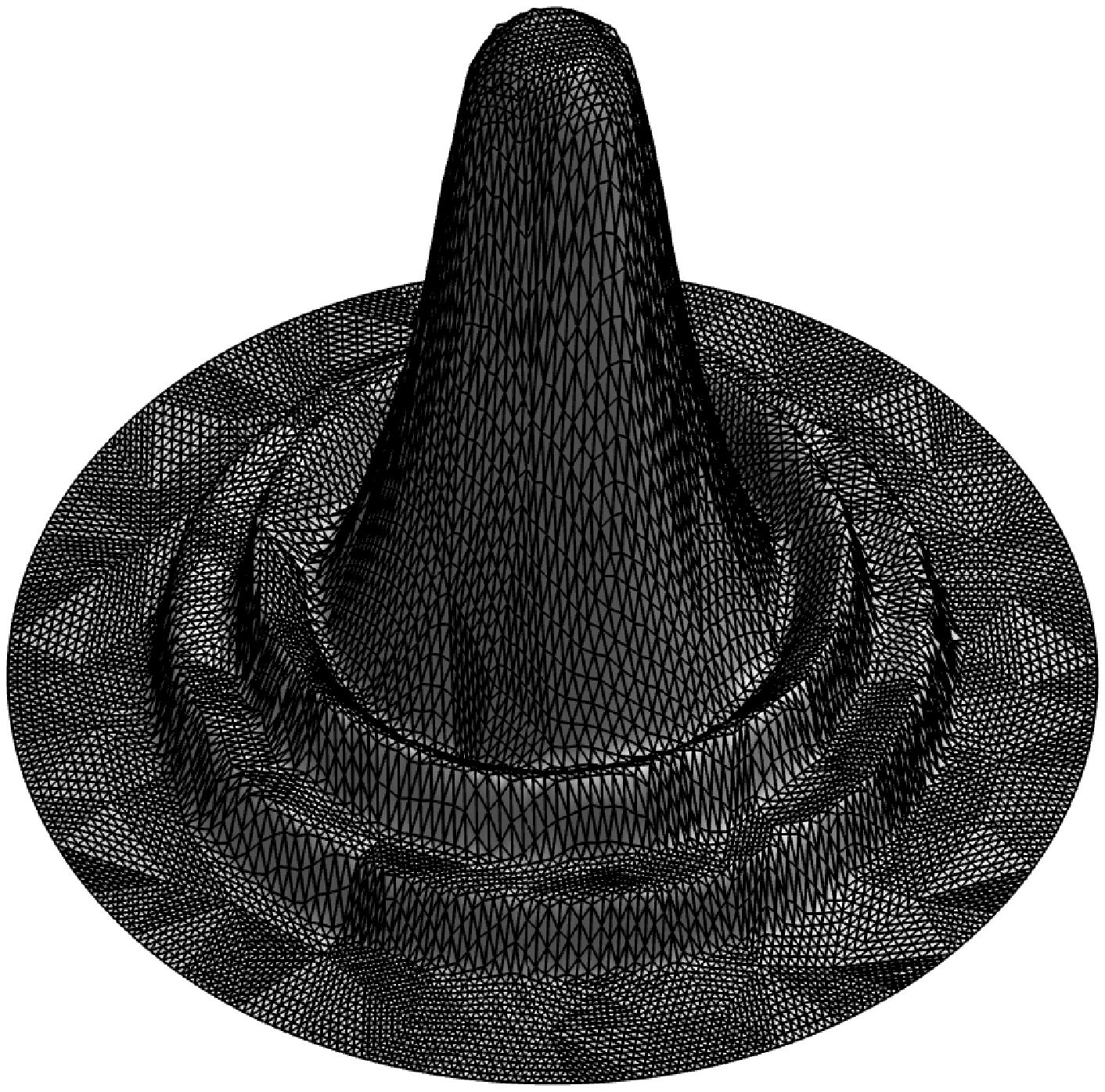}
\par\end{center}

\begin{center}
\protect\caption{\label{fig:Explosion_3D}Cross section of numerical solution in $xy$
plane of 3-D explosion test by third-order WLS-ENO at $t=0.1$. }

\par\end{center}%
\end{minipage}\hfill%
\begin{minipage}[t]{0.49\columnwidth}%
\begin{center}
\includegraphics[width=1\textwidth]{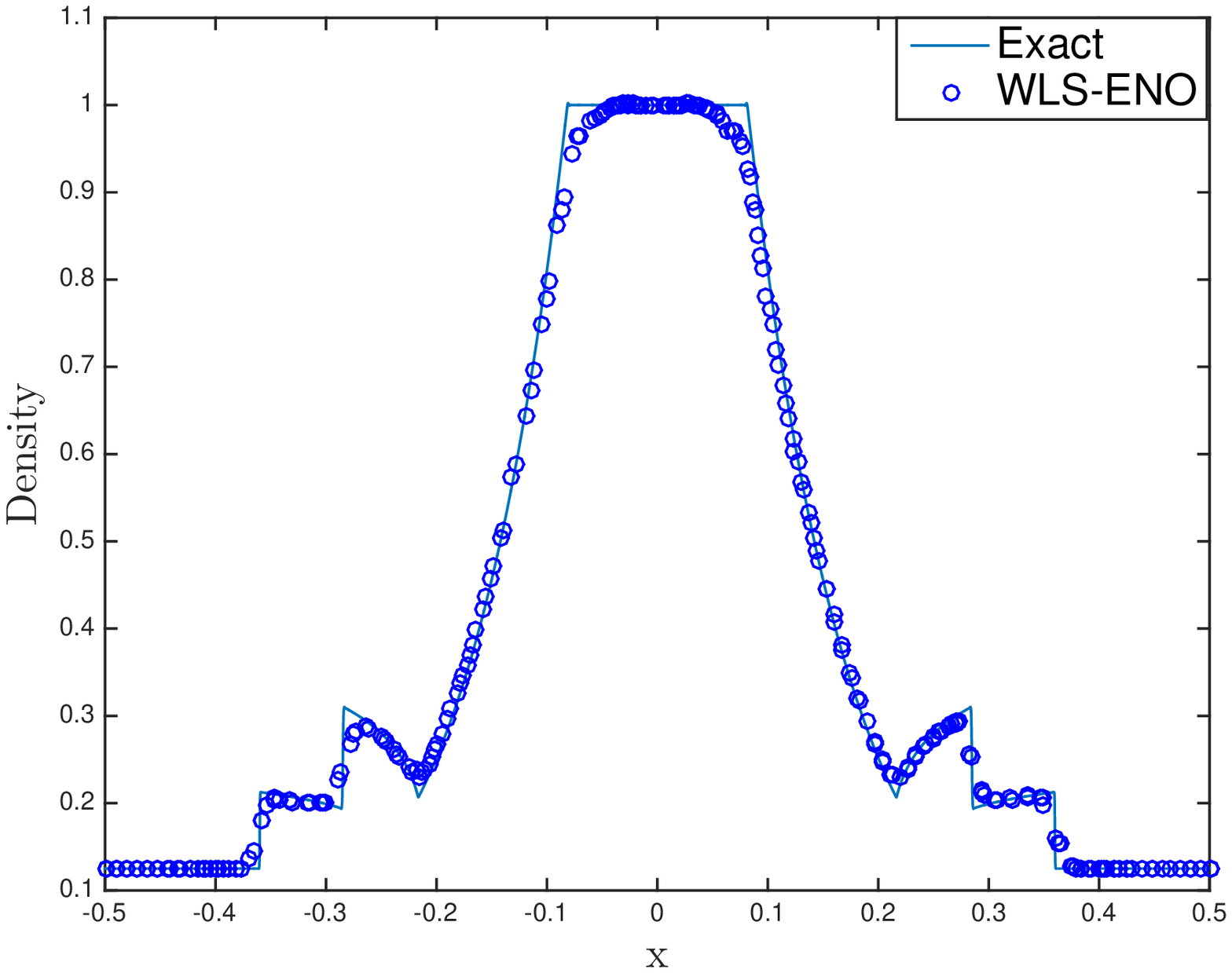}
\par\end{center}

\begin{center}
\protect\caption{\label{fig:1-D-comparision of Explosion 3D}Numerical solution along
$x$ axis of 3-D explosion test by third order WLS-ENO versus analytic
solution at $t=0.1$.}

\par\end{center}%
\end{minipage}
\end{figure}

We examine our scheme by solving an 3-D version of explosion test
problem. The initial condition of this problem is
\begin{equation}
\left(\rho,u,v,w,p\right)^{T}=\begin{cases}
\left(1,0,0,0,1\right)^{T} & \sqrt{x^{2}+y^{2}+z^{2}}\leq0.2\\
\left(0.125,0,0,0,0.1\right)^{T} & \sqrt{x^{2}+y^{2}+z^{2}}>0.2
\end{cases}.
\end{equation}
The computational domain is a unit ball centered at the origin, and
is tessellated by a tetrahedral mesh. We solve the problem up to $t=0.1$.
We computed the reference solution based on (\ref{eq:Explosion 1D})
for comparison with our numerical results. Figures~\ref{fig:Explosion_3D}
and \ref{fig:1-D-comparision of Explosion 3D} show the numerical
solutions, which agree with the analytic solutions very well. 

\section{\label{sec:Conclusions-and-Future}Conclusions and Future Work}

In this paper, we introduced a new family of essentially non-oscillatory
schemes, called \emph{WLS-ENO}, in the context of finite volume methods
for solving hyperbolic conservation laws. The schemes are derived
based on Taylor series expansions and solved with a weighted least
squares formulation. They can be applied to both structured and unstructured
meshes. Over structured meshes, we showed that WLS-ENO delivers similar
and even better accuracy compared to WENO, while enabling a larger
stability region. For unstructured meshes, we showed that WLS-ENO
enables accurate and stable solutions. Its accuracy and stability
are rooted in the facts that the convexity requirement is satisfied
automatically in WLS-ENO, and the stencil can be adapted more easily
to ensure the well-conditioning of the approximations. We presented
detailed analysis of WLS-ENO in terms of accuracy in 2-D, and its
stability in the context of hyperbolic conservation laws in 1-D. We
also assessed the WLS-ENO with a large collection of test problems
in 1-D, 2-D, and 3-D, including wave equations, Burgers' equation,
and the Euler equations with fairly complicated initial conditions.
Our numerical results demonstrated that WLS-ENO is accurate and stable
over unstructured meshes for very complex problems.

In its current form, WLS-ENO still has some limitations. Its primary
disadvantage is that it has higher computational cost compared to
the traditional WENO scheme over structured meshes. However, for engineering
applications involving complex geometries, WLS-ENO provides a more
general tool for dealing with piecewise smooth functions. To improve
efficiency, it is also possible to develop a hybrid method that utilizes
the traditional WENO on structured meshes in the interior and utilizes
WLS-ENO over unstructured meshes near complex boundaries. Another
limitation of the present formulation of WLS-ENO is that it only applies
to finite volume methods. As future work, we will extend WLS-ENO to
support generalized finite difference methods on unstructured meshes,
optimize the performance of WLS-ENO, and apply the methods to applications
in computational fluid dynamics.

\section*{Acknowledgements}

This work was supported by DoD-ARO under contract \#W911NF0910306.
The second author is also supported by a subcontract to Stony Brook
University from Argonne National Laboratory under Contract DE-AC02-06CH11357
for the SciDAC program funded by the Office of Science, Advanced Scientific
Computing Research of the U.S. Department of Energy.

\section*{References}

\bibliographystyle{elsarticle-num}
\bibliography{Refs/refs,Refs/WLS_ENO}

\end{document}